%% file: main.tex
\documentclass{siamart250211}

\usepackage[T1]{fontenc}    
\usepackage{lmodern}        
\usepackage{subcaption}     
\usepackage{siunitx}        
\usepackage{xfrac}          
\usepackage{booktabs}       

\crefname{subsection}{section}{sections} 

\usepackage{tikz}
\usetikzlibrary{arrows}     
\usetikzlibrary{positioning} 

\input{preamble/maths}      
\input{preamble/text}       

\begin{document}

\input{sections/frontmatter}

\input{sections/introduction}

\input{sections/model-definition}
\input{sections/pullout-bar}
\input{sections/three-point-hole}
\input{sections/conclusion}

\input{sections/backmatter}

\end{document}

%% file: preamble/maths.tex
\input{preamble/preamble/maths}

\newcommand*{\domain}{\Omega}
\newcommand*{\boundary}{\partial\Omega}
\newcommand*{\dirichletboundary}{\Gamma_\text{D}}
\newcommand*{\neumannboundary}{\Gamma_\text{N}}
\newcommand*{\differentialoperator}{\opD}
\newcommand*{\bilinearoperator}{\opA}

\newcommand*{\parameter}{\theta}
\newcommand*{\parameters}{\vtheta}
\newcommand*{\parametercount}{k}
\newcommand*{\parametrizeddifferentialoperator}{\differentialoperator_\parameters}
\newcommand*{\parametrizedbilinearoperator}{\bilinearoperator_\parameters}

\newcommand*{\hyperparameters}{\veta}

\newcommand*{\ofx}{{(x)}}
\newcommand*{\ofvx}{{(\vx)}}

\newcommand*{\ofxxprime}{{(x, x')}}
\newcommand*{\ofOmega}{{(\Omega)}}

\newcommand*{\solutionsymbol}{u}
\newcommand*{\sourcesymbol}{f}
\newcommand*{\solution}{\solutionsymbol}
\newcommand*{\source}{\sourcesymbol}
\newcommand*{\solutionvector}{\vec{\solutionsymbol}}

\newcommand*{\sourcevector}{\vec{\sourcesymbol}}

\newcommand*{\stiffnessmatrix}{\mK}

\newcommand*{\testfunctionsymbol}{\phi}

\newcommand*{\testfunction}{\testfunctionsymbol}

\newcommand*{\testvector}{\vphi}

\newcommand*{\testfunctioncount}{n}

\newcommand*{\elementcount}{n_e}

\newcommand*{\arbitraryfunctionv}{v}
\newcommand*{\arbitraryfunctionw}{w}
\newcommand*{\arbitraryvectorv}{\vv}
\newcommand*{\arbitraryvectorw}{\vw}

\NewDocumentCommand{\Lspace}{O{2}}{L^{#1}}
\NewDocumentCommand{\Hspace}{s m}{\IfBooleanTF{#1}{H^{#2}_0}{H^{#2}}}

\newcommand*{\reals}[1][]{\ifstrempty{#1}{\sR}{\sR^{#1}}}
\newcommand*{\dimensionality}{d}

\newcommand*{\meshsymbol}{\mathcal{T}}
\newcommand*{\meshsize}{h}
\newcommand*{\fem}[1]{{#1}^\text{\meshsize}}

\newcommand*{\femsolution}{\fem{\solution}}
\newcommand*{\mesh}{\fem{\meshsymbol}}
\newcommand*{\solutionspace}{\opV}
\newcommand*{\testspace}{\solutionspace_0}
\newcommand*{\trialspace}{\solutionspace}
\newcommand*{\femsolutionspace}{\fem{\solutionspace}}
\newcommand*{\femtestspace}{\fem{\solutionspace}_0}

\newcommand*{\weakformoperator}{\opL}

\newcommand*{\displacementfield}{\vu}
\newcommand*{\homdisplacementfield}{\vu_0}
\newcommand*{\femhomdisplacementfield}{\fem{\vu}_0}
\newcommand*{\inhomdisplacementfield}{\hat{\vu}}
\newcommand*{\strainfield}{\vvarepsilon}
\newcommand*{\stressfield}{\vsigma}
\newcommand*{\sourcefield}{\vf}
\newcommand*{\tractionfield}{\vt}
\newcommand*{\normalfield}{\vn}
\newcommand*{\elasticitytensor}{C}
\newcommand*{\dirichletfield}{\displacementfield_\text{D}}
\newcommand*{\neumannfield}{\tractionfield_\text{N}}

\newcommand*{\gradientof}{\nabla}
\newcommand*{\divergenceof}{\nabla \cdot}
\newcommand*{\symmetricgradientof}{\gradientof_\text{s}}

\newcommand*{\observations}{\vy}
\newcommand*{\observationcount}{m}
\newcommand*{\observationoperator}{\opP}
\newcommand*{\observationlocations}{\mX}
\newcommand*{\prior}{\probabilityof{\parameters}}
\newcommand*{\likelihood}{\probabilityof[\parameters]{\observations}}
\newcommand*{\posterior}{\probabilityof[\observations]{\parameters}}
\newcommand*{\evidence}{\probabilityof{\observations}}

\newcommand*{\noise}{\ve}
\newcommand*{\noisescale}{\sigma_e}
\newcommand*{\identitymatrix}{\mI}
\newcommand*{\noisecovariance}{\noisescale^2 \identitymatrix}
\newcommand*{\significancelevel}{\alpha}

\newcommand*{\bfemprior}{\probabilityof[\parameters]{\solution}}
\newcommand*{\bfemposterior}{\probabilityof[\sourcevector,\parameters]{\solution}}

\newcommand*{\bfemscale}{\sigma_u}
\newcommand*{\argmaxbfemscale}{\hat{\sigma}_u}

\newcommand*{\bfemposteriormean}{m_*}
\newcommand*{\bfemposteriorcovariance}{\opS_*}

\newcommand*{\covariancefunction}{k}

\newcommand*{\perturbed}[1]{\tilde{#1}}
\newcommand*{\nodelocation}{\vx}
\newcommand*{\rmfemscale}{p}
\newcommand*{\perturbation}{\valpha}
\newcommand*{\perturbations}{\valpha}
\newcommand*{\perturbedmesh}{\fem{\perturbed{\meshsymbol}}}
\newcommand*{\perturbedfemsolution}{\fem{\perturbed{\solutionsymbol}}}

\newcommand*{\misspecification}{\vd}
\newcommand*{\misspecificationkernel}{k_d}
\newcommand*{\misspecificationlength}{\ell_d}
\newcommand*{\misspecificationscale}{\sigma_d}
\newcommand*{\statfemscale}{\rho}

\newcommand*{\samplecount}{N}
\newcommand*{\burninsamplecount}{\samplecount_\text{burn}}
\newcommand*{\pseudomarginalsamplecount}{M}
\newcommand*{\tempering}{\tau}

\newcommand*{\youngsmodulus}{E}
\newcommand*{\crosssection}{A}
\newcommand*{\barstiffness}{E\mspace{-2mu}A}
\newcommand*{\springstiffness}{k}
\newcommand*{\pulloutforce}{F}
\newcommand*{\wavenumber}{\nu}
\newcommand*{\poissonsratio}{\nu}

\newcommand*{\beamheight}{H}
\newcommand*{\beamlength}{L}
\newcommand*{\beamoverhang}{c}
\newcommand*{\holelocx}{x}
\newcommand*{\holelocy}{y}
\newcommand*{\holesize}{d}
\newcommand*{\holeangle}{\alpha}
\newcommand*{\holeradius}{r}
\newcommand*{\observationspacing}{h_\text{m}}
\newcommand*{\prescribeddisplacement}{w}

%% file: preamble/preamble/maths.tex
\usepackage{amsmath}        
\usepackage{amsfonts}       
\usepackage{bm}             
\usepackage{mathcommand}    
\usepackage{mathtools}      
\usepackage{xparse}         

\renewcommand*{\vec}[1]{{\bm{#1}}}

\newcommand*{\vd}{\vec{d}}
\newcommand*{\ve}{\vec{e}}
\newcommand*{\vf}{\vec{f}}

\newcommand*{\vn}{\vec{n}}

\newcommand*{\vt}{\vec{t}}
\newcommand*{\vu}{\vec{u}}
\newcommand*{\vv}{\vec{v}}
\newcommand*{\vw}{\vec{w}}
\newcommand*{\vx}{\vec{x}}
\newcommand*{\vy}{\vec{y}}

\newcommand*{\valpha}{\vec{\alpha}}

\newcommand*{\vvarepsilon}{\vec{\varepsilon}}

\newcommand*{\veta}{\vec{\eta}}
\newcommand*{\vtheta}{\vec{\theta}}

\newcommand*{\vsigma}{\vec{\sigma}}

\newcommand*{\vphi}{\vec{\phi}}

\newcommand*{\vnull}{\vec{0}}

\newcommand*{\mat}[1]{{\bm{#1}}}

\newcommand*{\mI}{\mat{I}}

\newcommand*{\mK}{\mat{K}}

\newcommand*{\mX}{\mat{X}}

\newcommand*{\op}[1]{{\mathcal{#1}}}

\newcommand*{\opA}{\op{A}}

\newcommand*{\opD}{\op{D}}

\newcommand*{\opL}{\op{L}}

\newcommand*{\opP}{\op{P}}

\newcommand*{\opS}{\op{S}}

\newcommand*{\opV}{\op{V}}

\newcommand*{\set}[1]{{\mathbb{#1}}}

\newcommand*{\sR}{\set{R}}

\newcommand*{\super}[2]{{#1^{#2}}}

\newcommand*{\transposesymbol}{\mathsf{T}}
\newcommand*{\hermitiansymbol}{\mathsf{H}}
\newcommand*{\adjointsymbol}{*}
\newcommand*{\inversesymbol}{-}

\newcommand*{\transpose}[1]{\super{#1}{\transposesymbol}}

\newcommand*{\hermitian}[1]{\super{#1}{\hermitiansymbol}}
\newcommand*{\adjoint}[1]{\super{#1}{\adjointsymbol}}
\newcommand*{\inverse}[1]{\super{#1}{\inversesymbol 1}}

\newcommand*{\T}[1]{\transpose{#1}}
\renewmathcommand{\H}[1]{\hermitian{#1}}
\newcommand*{\adj}[1]{\adjoint{#1}}
\newcommand*{\inv}[1]{\inverse{#1}}

\newcommand*{\Tinv}[1]{\super{#1}{\inversesymbol\transposesymbol}}
\newcommand*{\adjinv}[1]{\super{#1}{\inversesymbol\adjointsymbol}}

\newcommand*{\invT}[1]{\Tinv}
\newcommand*{\invadj}[1]{\adjinv}
\newcommand*{\pinvT}[1]{\Tpseudoinv}
\newcommand*{\pinvadj}[1]{\adjpseudoinv}

\DeclarePairedDelimiter{\parens}{\lparen}{\rparen}
\DeclarePairedDelimiter{\bracks}{\lbrack}{\rbrack}
\DeclarePairedDelimiter{\braces}{\lbrace}{\rbrace}
\DeclarePairedDelimiter{\angles}{\langle}{\rangle}

\DeclarePairedDelimiter{\doublepipes}{\|}{\|}

\newcommand{\bsmallmat}[1]{
\begin{bsmallmatrix}
    #1
\end{bsmallmatrix}
}

\newcommand*{\dif}[2][1]{\mathrm{d}\ifstrequal{#1}{1}{}{^{#1}} #2}
\newcommand*{\pdif}[2][1]{\partial\ifstrequal{#1}{1}{}{^{#1}} #2}
\newcommand*{\diff}[2][1]{\mathrm{d} #2\ifstrequal{#1}{1}{}{^{#1}}}
\newcommand*{\pdiff}[2][1]{\partial #2\ifstrequal{#1}{1}{}{^{#1}}}

\newcommand*{\dx}{\diff{x}}

\newcommand*{\integrateover}[3][]{\ifstrempty{#1}{\int}{\int_{#1}} {#2} \, \diff{#3}}

\NewDocumentCommand{\uniform}{s m m}{
  \IfBooleanTF{#1}
    {\mathcal{U}\parens*{#2, #3}}
    {\mathcal{U}\parens{#2, #3}}
}

\NewDocumentCommand{\normal}{s m m}{
  \IfBooleanTF{#1}
    {\mathcal{N}\parens*{#2, #3}}
    {\mathcal{N}\parens{#2, #3}}
}

\NewDocumentCommand{\dirac}{s m}{
  \IfBooleanTF{#1}
    {\delta\parens*{#2}}
    {\delta\parens{#2}}
}

\NewDocumentCommand{\GP}{s m m}{
  \IfBooleanTF{#1}
    {\mathcal{GP}\parens*{#2, #3}}
    {\mathcal{GP}\parens{#2, #3}}
}

\DeclareMathOperator{\variance}{Var}
\DeclareMathOperator{\covariance}{Cov}

\NewDocumentCommand{\probabilityof}{s o m}{
  \IfBooleanTF{#1}
    {p\parens*{#3\IfValueT{#2}{|#2}}}
    {p\parens{#3\IfValueT{#2}{|#2}}}
}

\NewDocumentCommand{\proposalof}{s m}{
  \IfBooleanTF{#1}
    {q\parens*{#2}}
    {q\parens{#2}}
}

\NewDocumentCommand{\expectationof}{s m}{
  \IfBooleanTF{#1}
    {\mathbb{E}\bracks*{#2}}
    {\mathbb{E}\bracks{#2}}
}

\NewDocumentCommand{\varianceof}{s m}{
  \IfBooleanTF{#1}
    {\variance\parens*{#2}}
    {\variance\parens{#2}}
}

\NewDocumentCommand{\covarianceof}{s m m}{
  \IfBooleanTF{#1}
    {\covariance\parens*{#2, #3}}
    {\covariance\parens{#2, #3}}
}

\DeclareMathOperator*{\diag}{diag}
\DeclareMathOperator*{\argmax}{arg\,max}
\DeclareMathOperator*{\argmin}{arg\,min}
\DeclareMathOperator*{\linspan}{span}
\DeclareMathOperator*{\rank}{rank}

\NewDocumentCommand{\diagof}{s m}{
  \IfBooleanTF{#1}
    {\diag\parens*{#2}}
    {\diag\parens{#2}}
}

\NewDocumentCommand{\argmaxof}{s o m}{
  \IfBooleanTF{#1}
    {\argmax\IfValueT{#2}{_{#2}}\parens*{#3}}
    {\argmax\IfValueT{#2}{_{#2}}\parens{#3}}
}

\NewDocumentCommand{\argminof}{s o m}{
  \IfBooleanTF{#1}
    {\argmin\IfValueT{#2}{_{#2}}\parens*{#3}}
    {\argmin\IfValueT{#2}{_{#2}}\parens{#3}}
}

\NewDocumentCommand{\linspanof}{s m}{
  \IfBooleanTF{#1}
    {\linspan\parens*{#2}}
    {\linspan\parens{#2}}
}

\NewDocumentCommand{\rankof}{s m}{
  \IfBooleanTF{#1}
    {\rank\parens*{#2}}
    {\rank\parens{#2}}
}

\NewDocumentCommand{\norm}{s o m}{
  \IfBooleanTF{#1}
    {\doublepipes*{#3}\IfValueT{#2}{_{#2}}}
    {\doublepipes{#3}\IfValueT{#2}{_{#2}}}
}

\NewDocumentCommand{\innerproduct}{s o m m}{
  \IfBooleanTF{#1}
    {\angles*{#3, #4}\IfValueT{#2}{_{#2}}}
    {\angles{#3, #4}\IfValueT{#2}{_{#2}}}
}

\NewDocumentCommand{\collection}{s m m m}{
  \IfBooleanTF{#1}
    {\braces*{#2, #3, \dots, #4}}
    {\braces{#2, #3, \dots, #4}}
}

%% file: preamble/text.tex
\input{preamble/preamble/text}

\usepackage{xcolor}
\usepackage{pifont}

\definecolor{darkgreen}{rgb}{0.0, 0.7, 0.0}
\definecolor{darkred}{rgb}{0.7, 0.0, 0.0}

\newcommand*{\goodcheck}{\textcolor{darkgreen}{\ding{51}}}
\newcommand*{\badcheck}{\textcolor{darkred}{\ding{55}}}

%% file: preamble/preamble/text.tex
\usepackage{xspace}

\makeatletter
\newcommand{\etal}{\@ifstar\etalstar\etalnostar}
\newcommand{\etc}{\@ifstar\etcstar\etcnostar}
\newcommand{\iid}{\@ifstar\iidstar\iidnostar}
\makeatother

\newcommand{\etalstar}{et~al.}
\newcommand{\etalnostar}{et~al.\@\xspace}
\newcommand{\ie}{i.e.\@\xspace}
\newcommand{\eg}{e.g.\@\xspace}
\newcommand{\etcstar}{etc.}
\newcommand{\etcnostar}{etc.\@\xspace}
\newcommand{\iidstar}{i.i.d.}
\newcommand{\iidnostar}{i.i.d.\@\xspace}

%% file: sections/frontmatter.tex
\title{The Bayesian Finite Element Method in Inverse Problems: a Critical Comparison between Probabilistic Models for Discretization Error}

\author{Anne Poot\thanks{Faculty of Civil Engineering and Geosciences, Delft University of Technology (\email{a.poot-1@tudelft.nl}, \email{i.rocha@tudelft.nl}, \email{f.p.vandermeer@tudelft.nl}).}
\and Iuri Rocha\footnotemark[1]
\and Pierre Kerfriden\thanks{Centres des Mat\'eriaux, Mines Paris---PSL (\email{pierre.kerfriden@minesparis.psl.eu}).}
\and Frans van der Meer\footnotemark[1]}

\maketitle

\begin{abstract}
    When using the finite element method (FEM) in inverse problems, its discretization error can produce parameter estimates that are inaccurate and overconfident.
    The Bayesian finite element method (BFEM) provides a probabilistic model for the epistemic uncertainty due to discretization error.
    In this work, we apply BFEM to various inverse problems, and compare its performance to the random mesh finite element method (RM-FEM) and the statistical finite element method (statFEM), which serve as a frequentist and inference-based counterpart to BFEM.
    We find that by propagating this uncertainty to the posterior, BFEM can produce more accurate parameter estimates and prevent overconfidence, compared to FEM.
    Because the BFEM covariance operator is designed to leave uncertainty only in the appropriate space, orthogonal to the FEM basis, BFEM is able to outperform RM-FEM, which does not have such a structure to its covariance.
    Although inferring the discretization error via a model misspecification component is possible as well, as is done in statFEM, the feasibility of such an approach is contingent on the availability of sufficient data.
    We find that the BFEM is the most robust way to consistently propagate uncertainty due to discretization error to the posterior of a Bayesian inverse problem.
\end{abstract}








%% file: sections/introduction.tex
\section{Introduction}
\label{sec:introduction}
Inverse problems are abundant in science and engineering, ranging from climate modeling to medical imaging and from subsurface characterization to structural health monitoring.
Because they are typically ill-posed, solving an inverse problem requires many evaluations of the forward problem, which may involve the numerical solution of a partial differential equation (PDE).
The associated computational cost is severe, and the issue of limiting the cost can be targeted from various angles including approximation of the posterior \cite{blei_variational_2017,schillings_convergence_2020}, efficient sampling strategies \cite{hoffman_no-u-turn_2014,Bierkens2018}, model order reduction \cite{arridge_approximation_2006,cui_data-driven_2015} and machine learning techniques \cite{raissi_physics-informed_2019,li_fourier_2021}.
Recently, interest has increased in approaching numerical computation from a probabilistic point of view \cite{hennig_probabilistic_2015,cockayne_bayesian_2019}, which allows the epistemic uncertainty associated with limited computation to be propagated to the posterior.
In this work, we focus on probabilistic models for discretization error of the finite element method (FEM), in the context of Bayesian inverse problems in solid mechanics.

In recent work \cite{poot_bayesian_2024}, we proposed the Bayesian finite element method (BFEM), which takes a Bayesian approach to solving PDEs and thereby provides a probabilistic model of discretization error.
A Gaussian process prior is assumed over the solution space, and conditioned on observations implied by the FEM shape functions.
We can contextualize BFEM as part of a family of Bayesian methods of weighted residuals, which encompasses a broad range of PDE solvers \cite{pfortner_physics-informed_2023}.
Instead of conditioning on FEM shape functions, one can condition on point evaluations to recover a Bayesian collocation method \cite{cockayne_probabilistic_2017} or condition on volume integrals to recover a Bayesian finite volume method \cite{weiland_scaling_2024}.
Bayesian physics-informed neural networks \cite{yang_b-pinns_2021} fall in a similar category: their key difference with Bayesian collocation methods lies in the parametrization of the prior, for which an artificial neural network is used.

On the other hand, various frequentist approaches to modeling solver error have been proposed, which typically rely on perturbing a traditional PDE solver in some manner to induce a distribution over the solution space.
Conrad \etal \cite{conrad_statistical_2017} modeled uncertainty due to discretization error by adding locally supported noise to the FEM basis, whereas the random mesh finite element method (RM-FEM) by Abdulle and Garegnani \cite{abdulle_probabilistic_2021} models the error by perturbing the node locations to randomize the basis itself.
Note that not all methods fall neatly into these two categories.
For example, various probabilistic solvers for ordinary differential equations (ODEs) based on Gaussian filtering techniques have been proposed \cite{kersting_active_2018,tronarp_probabilistic_2019}.
Although such approaches can be generalized to PDEs as well \cite{kramer_probabilistic_2022}, they cannot be applied to PDEs on arbitrary domains, which is why they are left out of the scope of this work.
For a broader overview and history of probabilistic solvers of ODEs and PDEs, the reader is referred to Hennig \etal \cite{hennig_probabilistic_2022}.

Finally, we remark that discretization error is one of many forms of model misspecification.
Several methods have been proposed that aim to capture model misspecification probabilistically.
Alberts and Bilionis \cite{alberts_physics-informed_2023} adapt information field theory \cite{enslin_information_2013} by adding a tempering parameter to the prior, which calibrates the level of trust in the physical modeling.
By learning this parameter, the reliability of the model can be inferred, and the uncertainty due to model misspecification can be propagated to quantities of interest.
It is also possible to account for model misspecification by adapting the likelihood itself, as is done in the statistical finite element method (statFEM) \cite{girolami_statistical_2021}.
Rather than assign a distribution directly to the solution field, a model misspecification component is added to the observation model, similar to \cite{kennedy_bayesian_2001}.
The uncertainty due to model misspecification (in the form of discretization error) can be modeled by learning the misspecification parameters.

The main contributions of this work are twofold:
first, it is the first application of BFEM to an inverse problem.
We verify that by modeling the discretization error probabilistically, and consistently propagating it to the posterior, the quality of the posterior can be improved.
Second, we compare BFEM against RM-FEM \cite{abdulle_probabilistic_2021} and statFEM \cite{girolami_statistical_2021}, which present alternatives to the Bayesian approach from the frequentist and model misspecification categories, respectively.
Despite the wide range of probabilistic PDE solvers that have been proposed over the years, few works exist to comparatively study their performance, and this work contains the first comparison of the these three methods in an inverse setting.

This paper is structured as follows:
in \cref{sec:model-definition}, we summarize the BFEM, RM-FEM and statFEM models and describe how they account for error by adapting the likelihood of a Bayesian inverse problem.
In \cref{sec:experiments}, we present two numerical experiments.
A pullout test is presented in \cref{subsec:pullout-bar} as a 1D, low-data, material inference problem, followed by a three-point bending test in \cref{subsec:three-point-hole} as a 2D, high-data, geometry inference problem.
Finally, we conclude the paper and discuss its implications in \cref{sec:conclusion}.

%% file: sections/model-definition.tex
\section{Model definition}
\label{sec:model-definition}
For clarity of exposition, we assume our forward problem to be given by a second-order PDE over a domain $\domain$ with homogeneous Dirichlet boundary conditions:
\begin{equation}\label{eq:strong-form}
    \begin{aligned}
        \parametrizeddifferentialoperator \solution\ofx &= \source\ofx & \text{ for } x &\in \domain \\
        \solution\ofx &= 0 & \text{ for } x &\in \boundary
    \end{aligned}
\end{equation}
Here, $\parametrizeddifferentialoperator$ is an elliptic, self-adjoint differential operator, which depends on a set of parameters $\parameters \in \reals[\parametercount]$.
Note that the problems of linear elasticity we are interested in are generally not described by \cref{eq:strong-form}, due to possibly inhomogeneous Dirichlet and Neumann boundary conditions and $\solution\ofx$ being a vector field in $\reals[\dimensionality]$.
Nonetheless, the core ideas of the methods presented remains the same, and we refer the reader to \cref{sec:elasticity} for a discussion on the specifics of linear elastic solid mechanics.

We now assume that we have some noisy measurements $\observations \in \reals[\observationcount]$ that depend on the solution of the PDE:
\begin{equation}\label{eq:observation-model}
    \observations = \observationoperator \solution\ofx + \noise
\end{equation}
where $\observationoperator: \Lspace\ofOmega \to \reals[\observationcount]$ is a linear observation operator that maps the solution to the observations, and $\ve \sim \normal*{0}{\noisecovariance}$.
We assume that $\observationoperator$ performs pointwise evaluation of $\solution$ at a set of observation points $\observationlocations \in \domain$, although generalizations to other linear observation operators are trivial.
Since $\solution\ofx$ is fully determined by $\parameters$, we can use this observation model directly to define the likelihood:
\begin{equation}\label{eq:likelihood}
    \likelihood = \normal*{\observationoperator \solution\ofx}{\noisecovariance}
\end{equation}
A prior distribution is assumed over the parameters $\prior$, after which Bayes' theorem gives us a posterior distribution $\posterior$ over the parameters:
\begin{equation}\label{eq:posterior}
    \posterior = \frac{\prior \likelihood}{\evidence}
\end{equation}
In \cref{subfig:graphical-models-exact}, a graphical model of the inference problem under exact numerics is shown.

\newcommand{\bayesiannet}[3]{
    \tikzset{
        latent/.style={circle, draw, minimum size=#3},
        obs/.style={circle, draw, fill=gray!30, minimum size=#3},
        const/.style={fill=black, circle, minimum size=3pt, inner sep=0pt}
    }

    \tikzset{
        every node/.style={anchor=center},
        every edge quotes/.style={anchor=center},
        on grid,
        node distance=#1, 
        baseline=(current bounding box.center),
        center of/.style={at=((1.center))}
    }

    \begin{figure}
        \centering
        \refstepcounter{figure}

        \begin{tikzpicture}
            \node at (0, -3.75*#1) {\refstepcounter{subfigure} (\alph{subfigure}) exact\label{subfig:graphical-models-exact}};
            \begin{scope}[shift={(0,0)}]
                \node[latent] (theta) at (0,-#1) {$\parameters$};
                \node[latent,right=of theta] (D) {$\parametrizeddifferentialoperator$};
                \node[latent,below=of D] (ux) {$\solution$};
                \node[const,left=of ux] (fx) {}; \node at (fx) [left=#2] {$\source$};
                \node[obs,below=of ux] (y) {$\observations$};
                \node[latent,left=of y] (e) {$\noise$};
                \node[const,left=of e] (se) {}; \node at (se) [left=#2] {$\noisescale$};

                \draw[->] (theta) -- (D);
                \draw[->] (D) -- (ux);
                \draw[->] (fx) -- (ux);
                \draw[->] (ux) -- (y);
                \draw[->] (e) -- (y);
                \draw[->] (se) -- (e);
            \end{scope}

            \node at (3.5*#1, -3.75*#1) {\refstepcounter{subfigure}(\alph{subfigure}) FEM\label{subfig:graphical-models-fem}};
            \begin{scope}[shift={(3.5*#1,0)}]
                \node[latent] (theta) {$\parameters$};
                \node[latent,right=of theta] (D) {$\parametrizeddifferentialoperator$};
                \node[latent,below=of D] (K) {$\stiffnessmatrix$};
                \node[const,left=of K] (Th) {}; \node at (Th) [above left=0cm] {$\mesh$};
                \node[latent,below=of K] (u) {$\solutionvector$};
                \node[const,left=of u] (f) {}; \node at (f) [above left=#2] {$\sourcevector$};
                \node[const,left=of f] (fx) {}; \node at (fx) [left=#2] {$\source$};
                \node[obs,below=of u] (y) {$\observations$};
                \node[latent,left=of y] (e) {$\noise$};
                \node[const,left=of e] (se) {}; \node at (se) [left=#2] {$\noisescale$};

                \draw[->] (theta) -- (D);
                \draw[->] (D) -- (K);
                \draw[->] (Th) -- (K);
                \draw[->] (K) -- (u);
                \draw[->] (Th) -- (f);
                \draw[->] (f) -- (u);
                \draw[->] (fx) -- (f);
                \draw[->] (u) -- (y);
                \draw[->] (e) -- (y);
                \draw[->] (se) -- (e);
            \end{scope}

            \node at (7*#1, -3.75*#1) {\refstepcounter{subfigure}(\alph{subfigure}) BFEM\label{subfig:graphical-models-bfem}};
            \begin{scope}[shift={(7*#1,0)}]
                \node[latent] (theta) {$\parameters$};
                \node[latent,below=of theta] (D) {$\parametrizeddifferentialoperator$};
                \node[latent,below=of D] (fx) {$\source$};
                \node[latent,right=of fx] (ux) {$\solution$};
                \node[obs,left=of fx] (f) {$\sourcevector$};
                \node[const,above=of f] (Th) {}; \node at (Th) [above=#2] {$\mesh$};
                \node[const,above=of ux] (su) {}; \node at (su) [above=#2] {$\bfemscale$};
                \node[obs,below=of ux] (y) {$\observations$};
                \node[latent,left=of y] (e) {$\noise$};
                \node[const,left=of e] (se) {}; \node at (se) [left=#2] {$\noisescale$};

                \draw[->] (theta) -- (D);
                \draw[->] (D) -- (ux);
                \draw[->] (D) -- (fx);
                \draw[->] (ux) -- (fx);
                \draw[->] (fx) -- (f);
                \draw[->] (Th) -- (f);
                \draw[->] (su) -- (ux);
                \draw[->] (ux) -- (y);
                \draw[->] (e) -- (y);
                \draw[->] (se) -- (e);
            \end{scope}

            \node at (1.25*#1, -4.75*#1) {\refstepcounter{subfigure}(\alph{subfigure}) RM-FEM\label{subfig:graphical-models-rmfem}};
            \begin{scope}[shift={(1.25*#1,-5.5*#1)}]
                \node[latent] (theta) {$\parameters$};
                \node[latent,right=of theta] (D) {$\parametrizeddifferentialoperator$};
                \node[latent,below=of D] (K) {$\stiffnessmatrix$};
                \node[latent,left=of K] (Th) {$\mesh$};
                \node[latent,left=of Th] (alpha) {$\perturbations$};
                \node[latent,below=of K] (u) {$\solutionvector$};
                \node[latent,left=of u] (f) {$\sourcevector$};
                \node[const,left=of f] (fx) {}; \node at (fx) [left=#2] {$\source$};
                \node[obs,below=of u] (y) {$\observations$};
                \node[latent,left=of y] (e) {$\noise$};
                \node[const,left=of e] (se) {}; \node at (se) [left=#2] {$\noisescale$};

                \draw[->] (theta) -- (D);
                \draw[->] (D) -- (K);
                \draw[->] (Th) -- (K);
                \draw[->] (alpha) -- (Th);
                \draw[->] (K) -- (u);
                \draw[->] (Th) -- (f);
                \draw[->] (f) -- (u);
                \draw[->] (fx) -- (f);
                \draw[->] (u) -- (y);
                \draw[->] (e) -- (y);
                \draw[->] (se) -- (e);
            \end{scope}

            \node at (5.25*#1, -4.75*#1) {\refstepcounter{subfigure}(\alph{subfigure}) statFEM\label{subfig:graphical-models-statfem}};
            \begin{scope}[shift={(4.75*#1,-5.5*#1)}]
                \node[latent] (theta) {$\parameters$};
                \node[latent,right=of theta] (D) {$\parametrizeddifferentialoperator$};
                \node[latent,below=of D] (K) {$\stiffnessmatrix$};
                \node[const,left=of K] (Th) {}; \node at (Th) [above left=0cm] {$\mesh$};
                \node[latent,below=of K] (u) {$\solutionvector$};
                \node[const,left=of u] (f) {}; \node at (f) [above left=#2] {$\sourcevector$};
                \node[const,left=of f] (fx) {}; \node at (fx) [left=#2] {$\source$};
                \node[obs,below=of u] (y) {$\observations$};
                \node[latent,left=of y] (e) {$\noise$};
                \node[const,left=of e] (se) {}; \node at (se) [left=#2] {$\noisescale$};
                \node[latent,right=of y] (d) {$\misspecification$};
                \node[latent,below=of y] (rho) {$\statfemscale$};
                \node[latent] (ld) at (d |- u) {$\misspecificationlength$};
                \node[latent] (sd) at (d |- rho) {$\misspecificationscale$};

                \draw[->] (theta) -- (D);
                \draw[->] (D) -- (K);
                \draw[->] (Th) -- (K);
                \draw[->] (K) -- (u);
                \draw[->] (Th) -- (f);
                \draw[->] (f) -- (u);
                \draw[->] (fx) -- (f);
                \draw[->] (u) -- (y);
                \draw[->] (e) -- (y);
                \draw[->] (se) -- (e);
                \draw[->] (d) -- (y);
                \draw[->] (ld) -- (d);
                \draw[->] (sd) -- (d);
                \draw[->] (rho) -- (y);
            \end{scope}
        \end{tikzpicture}
        \addtocounter{figure}{-1}

        \caption{
            Graphical models of all statistical models considered in this work.
            A circle is used to indicate a probabilistic variable.
            If this variable is observed, the circle is shaded.
            A dot is used to indicate a deterministic (input) variable.
            Conditional dependencies between variables are indicated by arrows.
            For more details on probabilistic graph models, the reader is referred to \cite[Chapter 8]{bishop_pattern_2006}.
        }
        \label{fig:graphical-models}
    \end{figure}
}

\bayesiannet{1.2cm}{0.1cm}{0.8cm}

\subsection{Finite element method}
\label{subec:fem-definition}
It is generally not possible to compute the exact solution $\solution\ofx$ given a differential operator $\parametrizeddifferentialoperator$ and a source term $\source\ofx$, and so we resort to FEM to numerically approximate $\solution\ofx$.
We multiply both sides of \cref{eq:strong-form} by an arbitrary function $\arbitraryfunctionv\ofx \in \solutionspace = \Hspace*{1}$ and perform integration by parts to obtain the weak form of the PDE:
\begin{equation}
    \label{eq:weak-form}
    \begin{aligned}
        \parametrizedbilinearoperator\parens{\solution\ofx, \arbitraryfunctionv\ofx} &= \integrateover[\domain]{\source\ofx \arbitraryfunctionv\ofx}{x} & \forall \arbitraryfunctionv \in \Hspace*{1}
    \end{aligned}
\end{equation}
where $\parametrizedbilinearoperator$ is the bilinear form associated with $\parametrizeddifferentialoperator$.
A mesh $\mesh$ with element size $\meshsize$ is assumed over the domain $\domain$, which implies a set of $\testfunctioncount$ locally supported shape functions $\testvector\ofx = \T{\bsmallmat{\testfunction_1 \ofx & \testfunction_2 \ofx & \dots & \testfunction_\testfunctioncount \ofx}}$, which span $\femsolutionspace \subset \solutionspace$.
This yields a discretized weak form:
\begin{equation}
    \label{eq:discretized-weak-form}
    \begin{aligned}
        \parametrizedbilinearoperator\parens{\solution\ofx, \testfunction_i\ofx} &= \integrateover[\domain]{\source\ofx \testfunction_i\ofx}{x} & \forall i \in 1, \dots, \testfunctioncount
    \end{aligned}
\end{equation}
The equation above has an infinite set of solutions, and so we search $\femsolution \in \femsolutionspace$ to arrive at the FEM solution.

Noting that we have only applied linear operations to arrive at \cref{eq:discretized-weak-form} from \cref{eq:strong-form}, we summarize these steps as a single linear operator $\weakformoperator:\solutionspace \to \reals[\testfunctioncount]$, which takes any function $\arbitraryfunctionv \in \solutionspace$ and integrates it against the shape functions.
Its adjoint $\adj{\weakformoperator}:\reals[\testfunctioncount] \to \femsolutionspace$ performs weighted summation of the shape functions by any vector $\arbitraryvectorw \in \reals[\testfunctioncount]$:
\begin{align}
    \weakformoperator \arbitraryfunctionv\ofx &= \int_{\domain} \testvector\ofx \arbitraryfunctionv\ofx \dx &
    \adj{\weakformoperator} \arbitraryvectorw &= \sum_{i=1}^\testfunctioncount w_i \testfunction_i\ofx
\end{align}
Using this weak form operator $\weakformoperator$, we can very succinctly denote the relationship between the differential operator $\parametrizeddifferentialoperator$ and stiffness matrix $\stiffnessmatrix$, as well as between the source term $\source\ofx$ and the source vector $\sourcevector$:
\begin{align}
    \stiffnessmatrix &= \weakformoperator \parametrizeddifferentialoperator \adj{\weakformoperator} &
    \sourcevector &= \weakformoperator \source\ofx
\end{align}
The FEM approximation $\femsolution\ofx$ can then be obtained by solving the linear system that arises and using the resulting solution vector $\solutionvector$ to weigh the shape functions:
\begin{align}
    \stiffnessmatrix \solutionvector &= \sourcevector &
    \femsolution\ofx &= \adj{\weakformoperator} \solutionvector
\end{align}
The likelihood under the FEM approximation is given by:
\begin{equation}\label{eq:likelihood-fem}
    \probabilityof*[\parameters]{\observations} = \normal*{\observationoperator \femsolution\ofx}{\noisecovariance}
\end{equation}

The graph model of the FEM approximation to the inverse problem is shown in \cref{subfig:graphical-models-fem}.
Comparing it to that of the exact inverse problem in \cref{subfig:graphical-models-exact}, we observe that the FEM approximation introduces an intermediate layer between $\parameters$ and $\observations$, which results in some loss of information.
As we will see, neglecting the effects of this approximation can lead to incorrect parameter estimates.
We will therefore present three potential methods to model this uncertainty, and propagate it to the posterior $\posterior$.

\subsection{Bayesian finite element method}
\label{subsec:bfem-definition}
The first probabilistic model for discretization error we consider is our recently-proposed BFEM \cite{poot_bayesian_2024}, which takes a Bayesian approach to modelling the FEM error.
A prior distribution is assumed over the solution field $\bfemprior$, whose covariance operator is chosen to be the inverse of the differential operator $\parametrizeddifferentialoperator$ (\ie, the integral of the Green's function of the PDE) with a scaling parameter $\bfemscale$.
With some abuse of notation\footnote{It is abuse of notation in the sense that the function $\solution$ itself does not follow a distribution, only each finite set of evaluations of the function $\solution(\mX)$ is jointly Gaussian.
For the purposes of this paper, this distinction is not crucial, and we will maintain this abuse of notation throughout the work.}, we write:
\begin{equation}\label{eq:bfem-prior}
    \probabilityof*[\parameters]{\solution} = \GP*{0}{\bfemscale^2 \inv{\parametrizeddifferentialoperator}}
\end{equation}
Note that, in contrast to the typical covariance function $\covariancefunction\ofxxprime$, which computes the covariance between two points in $\domain$ directly, the covariance operator takes two functions in $\Lspace\ofOmega$ and computes their covariance:
\begin{equation}\label{eq:bfem-covariance}
    \begin{aligned}
        \covarianceof{\arbitraryfunctionv}{\arbitraryfunctionw} &= \integrateover[\domain]{\bfemscale^2 \inv{\parametrizeddifferentialoperator} \arbitraryfunctionv\ofx \, \arbitraryfunctionw\ofx}{x} & \forall \arbitraryfunctionv, \arbitraryfunctionw &\in \Lspace(\domain)
    \end{aligned}
\end{equation}
We rewrite \cref{eq:discretized-weak-form} in terms of the weak form operator $\weakformoperator$:
\begin{equation}\label{eq:bfem-observation-model}
    \weakformoperator \parametrizeddifferentialoperator \solution\ofx = \weakformoperator \source\ofx
\end{equation}
and condition on it to obtain the BFEM posterior $\bfemposterior = \GP*{\bfemposteriormean}{\bfemposteriorcovariance}$ \cite{poot_bayesian_2024}, where
\begin{align}\label{eq:bfem-posterior-moments}
    \bfemposteriormean &= \adj{\weakformoperator} \inv{\parens{\weakformoperator \parametrizeddifferentialoperator \adj{\weakformoperator}}} \weakformoperator \source\ofx &
    \bfemposteriorcovariance &= \bfemscale^2 \parens{\inv{\parametrizeddifferentialoperator} - \adj{\weakformoperator} \inv{\parens{\weakformoperator \parametrizeddifferentialoperator \adj{\weakformoperator}}} \weakformoperator} \\
    &= \adj{\weakformoperator} \inv{\stiffnessmatrix} \sourcevector = \femsolution\ofx &
    &= \bfemscale^2 \parens{\inv{\parametrizeddifferentialoperator} - \adj{\weakformoperator} \inv{\stiffnessmatrix} \weakformoperator} \nonumber
\end{align}
For a more rigorous mathematical discussion of the well-definedness of conditioning a Gaussian process on an observation operator of this form, the reader is referred to the work of Pf\"ortner \etal \cite{pfortner_physics-informed_2023}.

Comparing the probabilistic graph model of BFEM in \cref{subfig:graphical-models-bfem} to that of FEM in \cref{subfig:graphical-models-fem}, two differences stand out.
First, we note that the edges connecting $\solution\ofx$ to $\source\ofx$ have switched direction.
Instead of assuming the source term $\source\ofx$ to be fixed, the prior that is assumed over the solution $\solution\ofx$ also implies a Gaussian process prior over $\source\ofx$ and a normal distribution over $\sourcevector$.
In the BFEM model, it is the combination of $\solution\ofx$ and $\parametrizeddifferentialoperator\ofx$ that determine $\source\ofx$, rather than the other way around.
Second, we note the absence of a stiffness matrix $\stiffnessmatrix$ in the graphical model, which at first glance appears contradictory to its appearance in \cref{eq:bfem-posterior-moments}.
However, the stiffness matrix only appears in the posterior moments under the specific choice of prior in \cref{eq:bfem-prior}.
Had we assumed another covariance function such as the squared exponential, the graph model would remain largely unchanged, but $\stiffnessmatrix$ would no longer have appeared in \cref{eq:bfem-posterior-moments}.
It is only this combination of prior and observation model that allows us to recover a natural connection to standard FEM in the expression for the posterior.

Since both $\bfemposterior$ and $\noise$ are Gaussian, the observation model in \cref{eq:observation-model} implies a Gaussian distribution over the observation vector $\observations$ as well:
\begin{equation}
    \probabilityof*[\sourcevector,\parameters]{\observations} = \normal*{\observationoperator \femsolution\ofx}{\bfemscale^2 \observationoperator \parens{\inv{\parametrizeddifferentialoperator} - \adj{\weakformoperator} \inv{\stiffnessmatrix} \weakformoperator} \adj{\observationoperator} + \noisecovariance}
\end{equation}
Note that under this choice of prior, the BFEM likelihood is the same as the FEM likelihood in \cref{eq:likelihood-fem}, except for an additional term in the covariance that adds uncertainty in the space orthogonal to the FEM basis.
The computational bottleneck of the BFEM approach lies in the evaluation of $\observationoperator \inv{\parametrizeddifferentialoperator} \adj{\observationoperator}$, which we will approximate using a hierarchically refined mesh.
We will use a single level of refinement throughout the work, meaning that 1D line elements get split in two and 2D triangle elements get split into four congruent triangle elements to generate a fine reference mesh.
For more practical details on implementing BFEM, the reader is referred to \cite{poot_bayesian_2024}.

The main purpose of the scaling hyperparameter $\bfemscale$ is to account for the magnitude of the source term in the uncertainty given by the posterior.
Since the PDE is assumed to be linear, a tenfold scaling of the source term implies a tenfold scaling of the solution and error as well.
To account for this, we infer $\bfemscale$ using the evidence approximation \cite[Section 3.5]{bishop_pattern_2006}:
for each evaluation of the forward problem, we choose $\bfemscale$ such that it maximizes the evidence $\probabilityof*[\parameters]{\sourcevector}$, given by:
\begin{equation}
    \probabilityof*[\parameters]{\sourcevector} = \normal*{\vnull}{\bfemscale^2 \stiffnessmatrix}
\end{equation}
Because the evidence is Gaussian, we can obtain a closed-form expression for its maximizer $\argmaxbfemscale$:
\begin{equation}
    \argmaxbfemscale = \argmaxof*[{\bfemscale \in \reals[+]}]{\probabilityof*[\parameters]{\sourcevector}} = \sqrt{\frac{\T{\sourcevector} \inv{\stiffnessmatrix} \sourcevector}{\testfunctioncount}} = \sqrt{\frac{\T{\sourcevector} \solutionvector}{\testfunctioncount}}
\end{equation}
Note that since $\bfemscale$ does not affect the posterior mean, we can cheaply compute it a posteriori using the energy norm of the FEM solution $\femsolution\ofx$.

This approach is much simpler than the rescaling of the eigenvalues we originally suggested in \cite[Section 3.4]{poot_bayesian_2024}.
Our motivation to opt for a simple scaling factor $\bfemscale$ here, rather than a full rescaling of the eigenvalues of the posterior covariance is threefold:
first, the computational cost of updating the eigendecomposition of the posterior covariance in every MCMC iteration is tremendous, and would make the numerical experiments in \cref{sec:experiments} infeasible.
Second, the rescaling is not necessarily consistent with the Bayesian framework, since there are no guarantees that there exists an equivalent rescaled prior to the rescaled posterior distribution.
Finally, the approach followed here allows us to contrast the effect of covariance structure on the posterior estimate to the effect of covariance magnitude.
The BFEM posterior $\bfemposterior$ has a very specific structure that spans the space orthogonal to the FEM solution.
By not modifying the BFEM posterior covariance, we test whether the presence of this structure benefits the inference procedure, compared to other models whose variance more directly resembles the error, but which do not have this covariance structure.

\subsection{Random mesh finite element method}
\label{subsec:rmfem-definition}
We now turn to RM-FEM \cite{abdulle_probabilistic_2021} as a frequentist contender to the Bayesian approach presented in the previous section.
Instead of relying on a single, fixed mesh $\mesh$, each node at location $\nodelocation_i$ is perturbed according to:
\begin{equation}
    \perturbed{\nodelocation}_i = \nodelocation_i + \meshsize_i^\rmfemscale \perturbation_i
\end{equation}
where $\rmfemscale$ is a scaling parameter that controls the magnitude of the perturbations, and $\perturbation_i$ is a random perturbation, sampled from a uniform distribution over an appropriately sized hyperball.
Following \cite[Theorem 2.9]{abdulle_probabilistic_2021}, we fix $\rmfemscale=1$ in this paper, which guarantees that the error between the unperturbed and perturbed FEM solutions $\norm{\perturbedfemsolution - \femsolution}$ converges at the same rate as the error of the FEM solution $\norm{\femsolution - \solution}$.
Boundary nodes are perturbed and projected back to the boundary of the domain, and nodes at observation locations are not perturbed at all, to avoid additional interpolation error, following our own recent finding on RM-FEM performance \cite{poot_effects_2025}.

In the probabilistic graphical model for RM-FEM in \cref{subfig:graphical-models-rmfem}, we see how the uncertainty introduced by the nodal perturbations $\perturbations$ propagates through the model:
each perturbation produces a different mesh $\perturbedmesh$, which via $\stiffnessmatrix$ and $\sourcevector$ produces a different FEM solution $\solutionvector$.
Note that, in contrast with the leaf nodes in the other models, we will not attempt to estimate $\perturbations$ from our observed data\footnote{This would be akin to selecting the mesh that most closely matches the data, and therefore minimizes the discretization error}, but rather leave $\perturbations$ as a latent variable in the model.
Therefore, to sample from the likelihood $\probabilityof[\parameters]{\observations}$, we need to marginalize over all possible mesh perturbations $\perturbations$:
\begin{equation}
    \probabilityof*[\parameters]{\observations} = \integrateover{\probabilityof*[\parameters,\perturbations]{\observations} \probabilityof*{\perturbations}}{\perturbations}
\end{equation}
Although evaluation of this integral is computationally intractable, two approaches to approximate it are presented by Garegnani \cite{garegnani_sampling_2021}.
The first, Monte Carlo within Metropolis, is a pseudomarginal method \cite{beaumont_estimation_2003,andrieu_pseudo-marginal_2009} which approximates the integral above by a sample estimate:
\begin{equation}\label{eq:monte-carlo-within-metropolis}
    \probabilityof*[\parameters]{\observations} \approx \frac{1}{\pseudomarginalsamplecount} \sum_{j=1}^{\pseudomarginalsamplecount} \probabilityof*[\parameters,\perturbations_j]{\observations}
\end{equation}
For each new proposal, $\pseudomarginalsamplecount$ new meshes are sampled independently for both the current and proposed state of $\parameters$ to recompute the acceptance ratio.
The pseudomarginal approach introduces some additional computational overhead, but has the benefit that all possible mesh perturbations are taken into account as $\samplecount \to \infty$.

The second approach, Metropolis within Monte Carlo, instead replaces the posterior by a sample approximation:
\begin{equation}\label{eq:metropolis-within-monte-carlo}
    \probabilityof*[\observations]{\parameters} \approx \frac{1}{\pseudomarginalsamplecount} \sum_{j=1}^{\pseudomarginalsamplecount} \probabilityof*[\observations,\perturbations_j]{\parameters}
\end{equation}
For both approaches, since each perturbed mesh can be generated independently, this sample approximation is trivially parallelizable.
In settings where the observation noise $\noisescale$ is small, \cref{eq:metropolis-within-monte-carlo} gives a rather rough, multimodal approximation of the true posterior.
For this reason, we will instead make use of the Monte Carlo within Metropolis approach from \cref{eq:monte-carlo-within-metropolis} in this work.

\subsection{Statistical finite element method}
\label{subsec:misspecification-definition}
Finally, we will consider statFEM \cite{girolami_statistical_2021} as a way to infer the discretization error from observations directly, rather than model it explicitly.
We note that it is not conventional in statFEM literature to consider discretization error to be part of the model misspecification.
More commonly, the FEM approximation of the statFEM posterior is assumed to be close enough to the true posterior, and this assumption is supported by convergence rate proofs \cite{papandreou_theoretical_2022,karvonen_error_2024}.
Nonetheless, we note that discretization error is not inherently different from any other form of model misspecification such as incorrectly assumed material properties, geometry, or loading conditions.
For this reason, we include statFEM in this work as an error-learning alternative to the error-modeling approaches that BFEM and RM-FEM represent.
The statFEM observation model is given by:
\begin{equation}\label{eq:statfem-observation-model}
    \observations = \statfemscale \observationoperator \solution\ofx + \misspecification + \noise
\end{equation}
Compared to \cref{eq:observation-model}, we have an additional scaling parameter $\statfemscale$ and a model misspecification component $\misspecification$, as initially proposed in \cite{kennedy_bayesian_2001}.
We will assume $\misspecification$ to be given by a zero-mean Gaussian process with a squared exponential covariance function, evaluated at the observation locations $\observationlocations$:
\begin{align}\label{eq:statfem-misspecification}
    \probabilityof*{\misspecification} &= \normal*{\vnull}{\misspecificationkernel(\observationlocations,\observationlocations)} &
    \misspecificationkernel(x,x') &= \misspecificationscale^2 \exp\parens*{-\frac{\norm{x-x'}^2}{2 \misspecificationlength^2}}
\end{align}

In \cref{subfig:graphical-models-statfem}, the graphical model used for statFEM is shown.
Note that, in contrast to BFEM and RM-FEM, no modifications to the model are introduced to any parent nodes of $\solutionvector$.
Rather, the difference between the statFEM and FEM graph models lies in the two additional parent nodes of $\observations$.
For the sake of consistency with the other models, we will maintain $\noisescale$ as a fixed, known parameter, and only learn $\statfemscale$, $\misspecificationlength$ and $\misspecificationscale$.
These three hyperparameters are gathered in the variable $\hyperparameters$, and learned them jointly with the parameters $\parameters$.
If we assume independence between $\parameters$ and $\hyperparameters$, the expression for the posterior in \cref{eq:posterior} becomes:
\begin{equation}
    \probabilityof*[\observations]{\parameters,\hyperparameters} = \frac{\probabilityof*{\parameters} \probabilityof*{\hyperparameters} \probabilityof*[\parameters, \hyperparameters]{\observations}}{\probabilityof*{\observations}}
\end{equation}

This is quite a simplistic model for discretization error, and various improvements are easily conceivable.
A broader set of prior distributions on the misspecification component $\misspecification$ could be explored like the Mat\'ern family \cite{rasmussen_gaussian_2005} and its SPDE approximations \cite{roininen_whittle-matern_2014,Lindgren2011}.
Alternatively, the discretization error can be taken into account via a separate additive component.
For example, in \cite{hermann_statistical_2024}, the statFEM model is enriched with an additional misspecification component, based on the full-field error estimate from \cite{rouse_probabilistic_2021}.
In a similar vein, one could combine the statFEM observation model with the BFEM error model, which would result in an additional BFEM term in the statFEM likelihood.
We limit the scope of this work to only a squared exponential covariance function, which is sufficient to study the difference in general behavior of an error inference in contrast with explicit error models like BFEM and RM-FEM.

\subsection{Model comparison}
It is useful to consider beforehand what constitutes good performance of a probabilistic model for discretization error, which is not immediately obvious.
We use three criteria to evaluate model performance.
First, we check for overconfidence of the model by considering a very coarse mesh for which the effect of discretization error on the posterior is significant.
We say that the model is overconfident if the ground truth falls outside the credible interval of the posterior distribution with significance level $\significancelevel$.
In this case, discretization error is not sufficiently accounted for by the model.
Second, we check for underconfidence of the model by considering a very fine mesh for which the effect of discretization error on the posterior is very small.
We say that the model is underconfident if a non-negligible difference between its posterior and the exact posterior remains, even if the discretization error itself is negligible.
In this case, the model accounts too much for discretization error, and the model is too uncertain about its prediction.
Third, we check for consistent convergence between these two extremes.
As the mesh is refined, the model should become more confident about its predictions and converge towards the exact posterior distribution.
We say that a model is performing well if it is not overconfident, not underconfident, and converges with mesh refinement.

\subsection{Random walk Metropolis sampling}
For all models and all examples in this work, we make use of the same random walk Metropolis \cite{metropolis_equation_1953} sampling scheme.
We start with a burn-in period of $\burninsamplecount=\num{10000}$ samples, during which we use a tempered likelihood:
\begin{equation}
    \posterior \propto \prior \likelihood^\tempering
\end{equation}
The tempering parameter $\tempering$ increases linearly from $\tempering=\num{0}$ at the first burn-in sample to $\tempering=\num{1}$ at the last.
We initialize the proposal distribution as a Gaussian whose mean and covariance equal that of the prior distribution, and adaptively scale its covariance during the burn-in period to maintain a good acceptance ratio.
This way, we ensure a well-calibrated proposal distribution and a reasonable starting point by the end of the burn-in period.
After the burn-in period, we discard all samples, fix the covariance of the proposal distribution and $\tempering=\num{1}$, and run the chain for an additional $\samplecount=\num{10000}$ samples.
For RM-FEM, we use the pseudomarginal approximation from \cref{eq:monte-carlo-within-metropolis} with $\pseudomarginalsamplecount=\num{100}$ in \cref{subsec:pullout-bar} and $\pseudomarginalsamplecount=\num{10}$ in \cref{subsec:three-point-hole}.

%% file: sections/pullout-bar.tex
\section{Numerical experiments}
\label{sec:experiments}
In this section, we compare the performance of BFEM, RM-FEM and statFEM using two different test cases.
Although both are academic examples, they are designed to resemble realistic applications within material sciences.
In \cref{subsec:pullout-bar}, we present a pullout test, where we aim to infer two material characteristics from only a single measurement.
At first glance, this may appear to be an unusually low-data setting; however, it is very common for material tests to only produce a load-displacement curve at the load point.
Typically, the entire curve is used to characterize the full nonlinear response of the material under loading, but for the elastic regime, only the initial slope can be used.
For a linear PDE, this is equivalent to conditioning on only a single force-displacement observation.

In \cref{subsec:three-point-hole}, we consider a three-point bending test, where we search for the location of a hole at some unknown location in the material.
This setup is a simplified version of digital image correlation methods that are used to detect material defects \cite{mccormick_digital_2010}.
Where digital image correlation infers defects in a 3D object by measuring (the difference between) displacements on its 2D boundary, we infer a defect in a 2D object by measuring displacements on its 1D boundary.
Although the parametrization of the defect as a single hole is not as realistic as, for example, crack formations at various unknown locations, it is more easily understood intuitively.
The code for both experiments is publically available at \href{https://gitlab.tudelft.nl/apoot1/probfem}{https://gitlab.tudelft.nl/apoot1/probfem}.

\subsection{Pullout test}
\label{subsec:pullout-bar}
\begin{figure}
    \centering
    \begin{subfigure}{0.5\textwidth}
        \centering
        \includegraphics[width=0.8\textwidth]{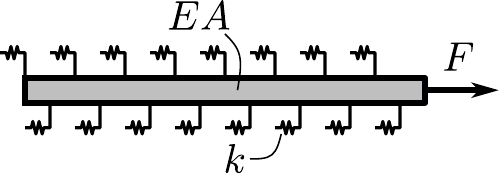}
        \caption{}
        \label{subfig:pullout-bar-overview}
    \end{subfigure}\hfill
    \begin{subfigure}{0.5\textwidth}
        \centering
        \includegraphics[width=\textwidth]{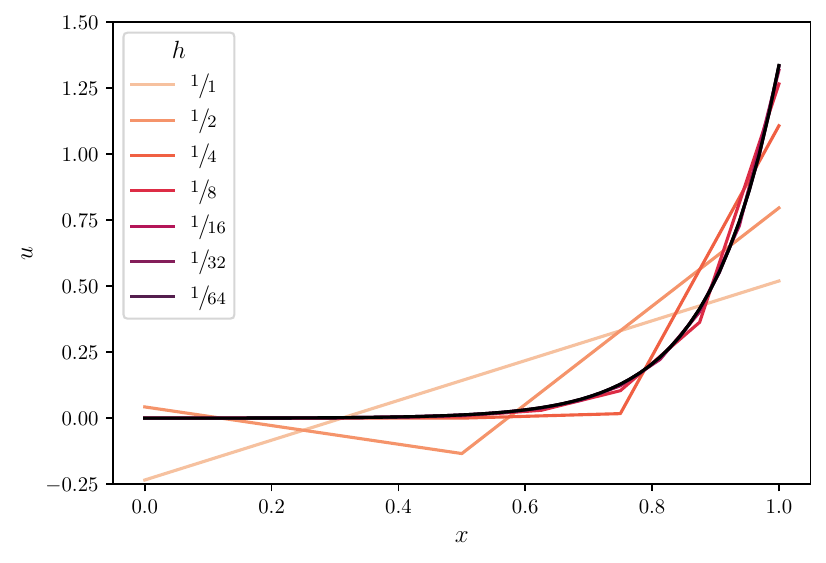}
        \caption{}
        \label{subfig:pullout-bar-exact}
    \end{subfigure}
    \caption{
        a)
        schematic overview of the pullout test.
        b)
        FEM solution $\femsolution\ofx$ for different mesh sizes $\meshsize$.
        The exact solution given by equation \cref{eq:pullout-bar-exact-solution} is shown in black.
    }
\end{figure}

We consider an elastic bar, embedded in an elastic medium.
The bar has a constant Young's modulus $\youngsmodulus$ and cross-sectional area $\crosssection$.
We model the medium as an elastic support along the bar with constant stiffness $\springstiffness$.
At the right end of the bar, a force $\pulloutforce$ is exerted.
In \cref{subfig:pullout-bar-overview}, a schematic overview of the problem is shown.
The pullout test is described by the following differential equation:
\begin{equation}
    \begin{aligned}
        -\frac{\diff{}}{\dx}\parens*{\barstiffness \, \frac{\diff{\solution}}{\dx}} + \springstiffness \, \solution &= 0 & \forall x &\in (0, 1) \\
        \frac{\diff{\solution}}{\dx}  &= 0 &\text{at } x &= 0 \\
        \frac{\diff{\solution}}{\dx}  &= \frac{\pulloutforce}{\barstiffness} &\text{at } x &= 1
    \end{aligned}
\end{equation}
which has a closed-form solution, given by:
\begin{equation}\label{eq:pullout-bar-exact-solution}
    \solution\ofx = \frac{\pulloutforce}{\sqrt{\springstiffness \, \barstiffness}} \frac{\cosh\parens*{\nu x}}{\sinh\parens*{\nu}}
\end{equation}
where $\wavenumber = \sqrt{\springstiffness/\barstiffness}$ is the inverse decay length of the system.
In \cref{subfig:pullout-bar-exact}, the exact solution $\solution\ofx$ is shown for $\barstiffness=\num{0.8}$, $\springstiffness=\num{70}$ and $\pulloutforce=\num{10}$ as a function of $x$.
The FEM solution $\femsolution\ofx$ is shown for various mesh densities, ranging from a single linear element to a uniform mesh consisting of 64 linear elements.

\begin{figure}
    \centering
    \begin{subfigure}{0.33\textwidth}
        \centering
        \includegraphics[width=\textwidth]{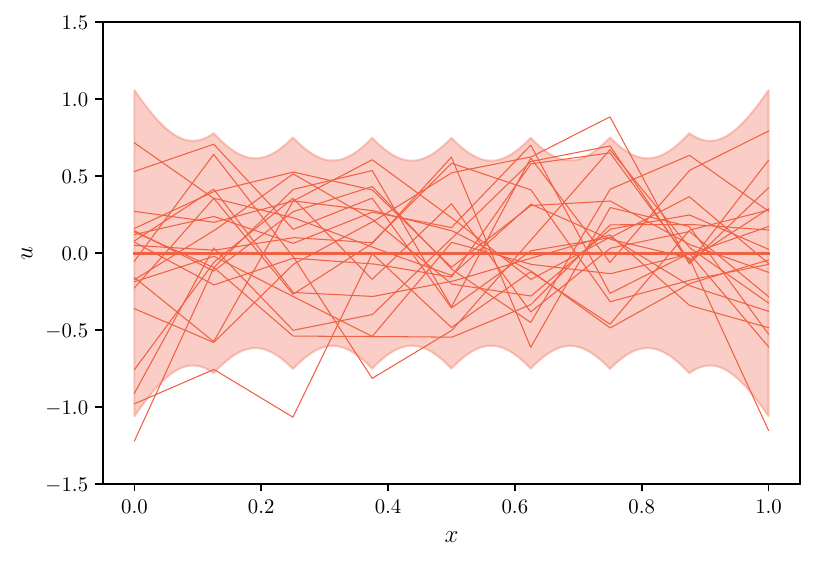}
        \caption{BFEM, $\probabilityof[\parameters]{\solution}$}
        \label{subfig:pullout-bar-bfem-prior-ref-8}
    \end{subfigure}\hfill
    \begin{subfigure}{0.33\textwidth}
        \centering
        \includegraphics[width=\textwidth]{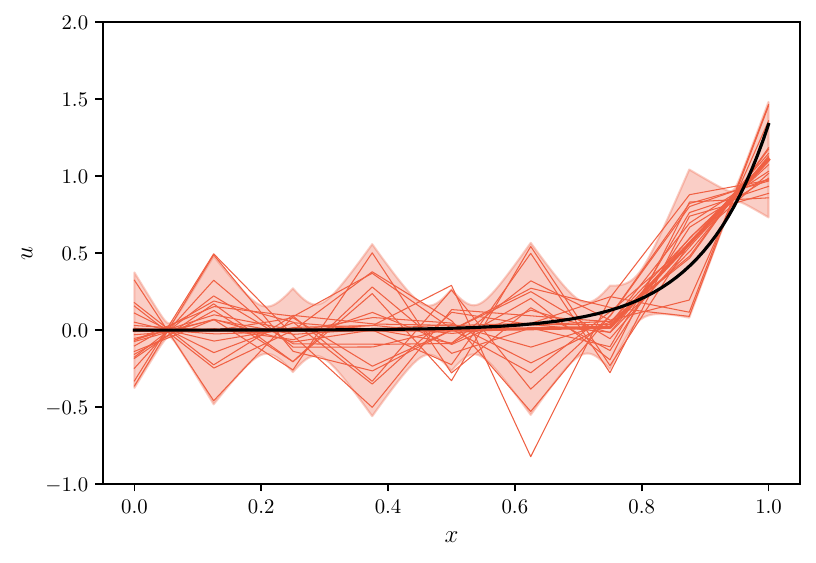}
        \caption{BFEM, $\probabilityof[\parameters,\sourcevector]{\solution}$}
        \label{subfig:pullout-bar-bfem-posterior-ref-8}
    \end{subfigure}\hfill
    \begin{subfigure}{0.33\textwidth}
        \centering
        \includegraphics[width=\textwidth]{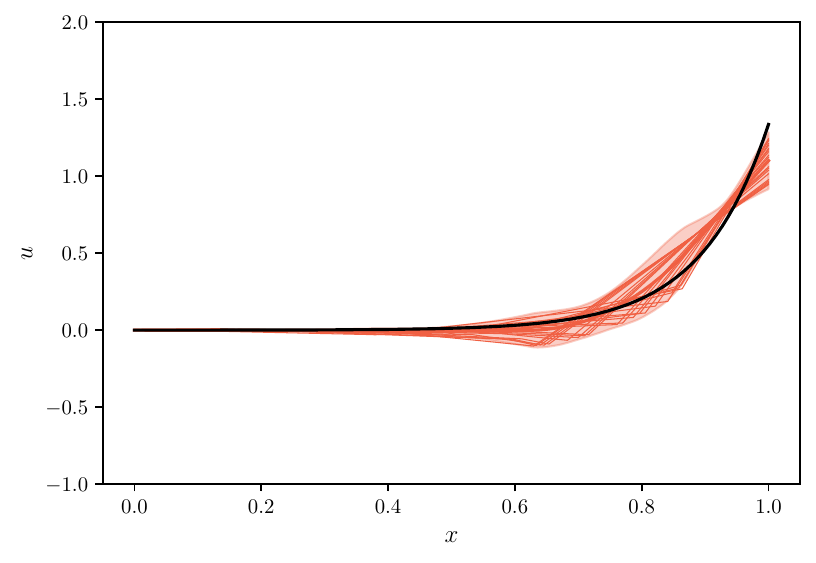}
        \caption{RM-FEM, $\probabilityof[\parameters]{\solution}$}
        \label{subfig:pullout-bar-rmfem-posterior}
    \end{subfigure}
    \caption{
        Distributions over the solution field $\solution\ofx$ for BFEM and RM-FEM.
    }
    \label{fig:pullout-bar-forward}
\end{figure}

We now consider what the distributions over $\solution\ofx$ look like for BFEM and RM-FEM\footnote{StatFEM is not included here, because it assumes a distribution over $\observations$ instead of $\solution\ofx$ (see \cref{fig:graphical-models}).}.
To investigate this, we consider the FEM solution where $\meshsize=\sfrac{1}{4}$, and plot the resulting BFEM and RM-FEM distributions in \cref{fig:pullout-bar-forward}.
For BFEM, we use a single level of refinement as a reference mesh (\ie, $\meshsize_f=\sfrac{1}{8}$), which yields the jagged pattern visible in the BFEM prior $\probabilityof[\parameters]{\solution}$ shown in \cref{subfig:pullout-bar-bfem-prior-ref-8}.
In the BFEM posterior $\probabilityof[\parameters,\sourcevector]{\solution}$, shown in \cref{subfig:pullout-bar-bfem-posterior-ref-8}, we see that this pattern persists, but subject to a contraction at the 5 observation nodes.
Note that, despite the FEM error being small to the left half of the bar, we still observe a large variance in this region.
This matches earlier observations we made in \cite[Section 3.4]{poot_bayesian_2024}, where we noted that the standard deviation of the BFEM posterior is not directly reflective of discretization error, because it only accounts for the structure of the FEM basis, not the specific loading conditions.

Turning to the RM-FEM results shown in \cref{subfig:pullout-bar-rmfem-posterior}, we see a distribution over $\solution\ofx$ that strongly differs from the one produced by BFEM.
The variance on the left half of the bar is very small, compared to the right half, which is in accordance with the small FEM error on the left half.
The contraction around $x=\sfrac{7}{8}$ is a consequence of the fact that the right-most node is a boundary node, and therefore not perturbed.
The variance at this node depends strongly on the location of the node to the left of it (at $x=\sfrac{3}{4}$), and their negative correlation causes a contraction between them.
The similarity between the FEM error and RM-FEM variance as a function of $x$ agrees with the fact that for RM-FEM, error estimators can be directly derived from the variance itself \cite[Section 3]{abdulle_probabilistic_2021}.
For BFEM, on the other hand, the full covariance structure needs to be taken into account (rather than only the diagonal) to retrieve the error from the posterior distribution \cite[Section 3.4]{poot_bayesian_2024}.

\subsubsection{Inverse problem}
\label{subsubsec:pullout-inverse}
\begin{figure}
    \centering
    \begin{subfigure}{0.5\textwidth}
        \centering
        \includegraphics[width=\textwidth]{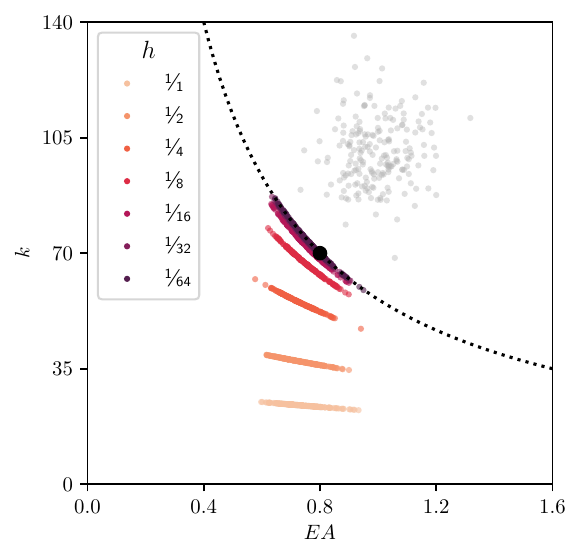}
        \caption{FEM}
        \label{subfig:pullout-bar-inverse-fem}
    \end{subfigure}\hfill
    \begin{subfigure}{0.5\textwidth}
        \centering
        \includegraphics[width=\textwidth]{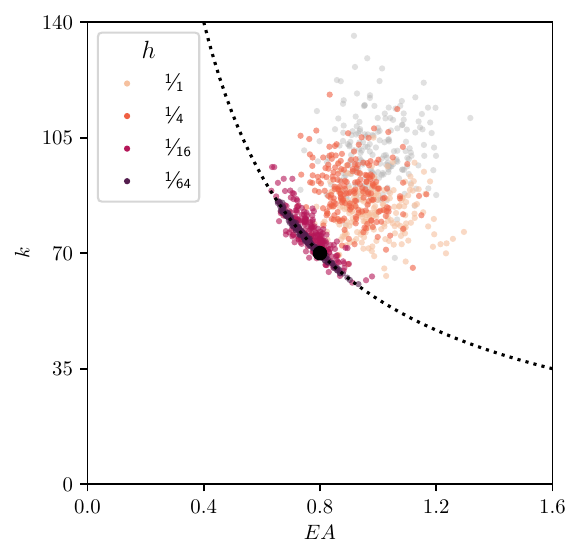}
        \caption{BFEM}
        \label{subfig:pullout-bar-inverse-bfem}
    \end{subfigure}

    \vspace{10pt}
    \begin{subfigure}{0.5\textwidth}
        \centering
        \includegraphics[width=\textwidth]{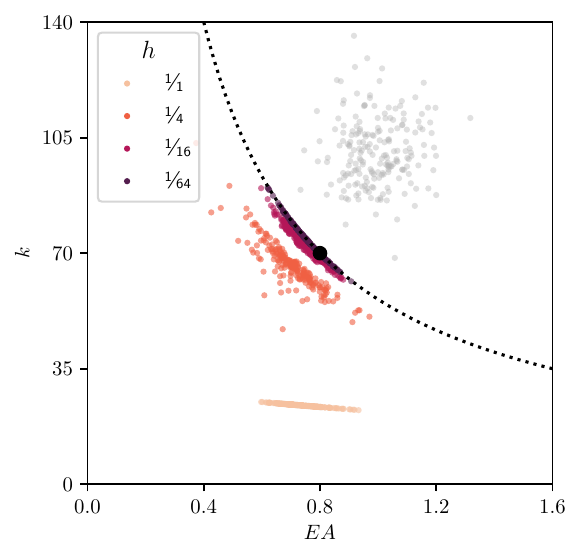}
        \caption{RM-FEM}
        \label{subfig:pullout-bar-inverse-rmfem}
    \end{subfigure}\hfill
    \begin{subfigure}{0.5\textwidth}
        \centering
        \includegraphics[width=\textwidth]{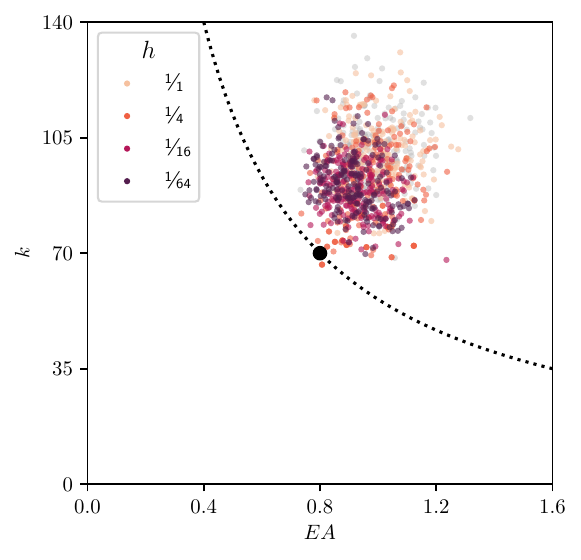}
        \caption{statFEM}
        \label{subfig:pullout-bar-inverse-statfem}
    \end{subfigure}
    \caption{
        Posterior distributions $\posterior$ for the pullout test inverse problem.
        For reference, samples from the prior distribution $\prior$ are shown in gray.
        The ground truth is indicated by a black point at the center of the plot.
        Along the black dotted line, $\springstiffness \barstiffness$ is constant.
    }
    \label{fig:pullout-bar-inverse}
\end{figure}

We now turn to an inverse setting, where we want to infer both $\barstiffness$ and $\springstiffness$.
To do so, we generate a single synthetic measurement by evaluating \cref{eq:pullout-bar-exact-solution} at $x=1$ for $\barstiffness=\num{0.8}$, $\springstiffness=\num{70}$ and $\pulloutforce=10$ and corrupting the resulting displacement with observation noise $\noisescale=10^{-3}$.
Clearly, the resulting inverse problem is ill-posed, as we aim to infer two parameters using only a single measurement.
This is where the benefit of using a Bayesian approach comes into play.
We assume the following prior distributions over $\barstiffness$ and $\springstiffness$:
\begin{align}
    \log \probabilityof*{\barstiffness} &= \normal*{\log 1}{0.1^2} &
    \log \probabilityof*{\springstiffness} &= \normal*{\log 100}{0.1^2}
\end{align}
This implies for $\barstiffness$ that \num{95.3}\% of the probability mass lies between \num{80} and \num{120}, and for $\springstiffness$ that \num{95.3}\% of the probability mass lies between \num{0.8} and \num{1.2}.
In \cref{fig:pullout-bar-inverse}, the resulting posterior distributions are shown for FEM, BFEM, RM-FEM and statFEM.

Looking at the FEM posterior distributions in \cref{subfig:pullout-bar-inverse-fem}, we see a clear convergence with mesh refinement to the exact posterior distribution.
The coarsest mesh $\meshsize=\sfrac{1}{1}$, which consists of just a single element, severely underestimates both $\barstiffness$ and $\springstiffness$.
This underestimation is a direct consequence of the underestimation of $\solution(x=1)$ shown in \cref{subfig:pullout-bar-exact} due to discretization error.
The FEM response is too stiff, which is compensated in the inverse setting by lowering the stiffness of the bar $\barstiffness$ and the stiffness of the springs $\springstiffness$.
As the mesh is refined, the error in the inferred parameters is reduced, until convergence is reached at $\meshsize\leq\sfrac{1}{32}$.

In \cref{subfig:pullout-bar-inverse-bfem}, we see that for the coarsest meshes, the BFEM posterior defaults to the prior.
As the mesh is refined, confidence about the FEM solution increases and the posterior gradually converges to the exact posterior distribution.
For the finest mesh, the FEM and BFEM posterior are virtually identical.
Note that for $\meshsize=\sfrac{1}{1}$, the reference mesh is the $\meshsize=\sfrac{1}{2}$ FEM mesh.
Although this reference mesh is not sufficiently fine to produce an accurate posterior (see \cref{subfig:pullout-bar-inverse-fem}), it is still able to inform us that the $\meshsize=\sfrac{1}{1}$ mesh cannot be trusted either.
This property could be useful in settings where computational budget is limited:
starting from a very coarse observation mesh, and a refined reference mesh, one can run the inverse problem.
If the resulting posterior is too wide, the refined mesh is used as an observation mesh, and refined to generate a new reference mesh.
This process can be repeated until the resulting posterior is sufficiently narrow, according to some stopping criterion.
This way, computational cost can be kept to a minimum while still remaining confident in the accuracy of the final posterior.

The RM-FEM posterior shown in \cref{subfig:pullout-bar-inverse-rmfem} strongly differs from the BFEM posterior.
Most notably, where BFEM converges from the prior towards the exact posterior with mesh refinement, RM-FEM converges from below similar to FEM.
For $\meshsize=\sfrac{1}{1}$, we recover exactly the FEM posterior, since all nodes are boundary nodes and therefore not perturbed.
For $\meshsize=\sfrac{1}{4}$, the RM-FEM posterior is closer to the ground truth than the FEM posterior and wider in the direction perpendicular to the manifold where $\springstiffness \barstiffness$ is constant.
We can understand this by returning to \cref{subfig:pullout-bar-rmfem-posterior}, and noting that the mesh perturbation produces a range of estimates of the solution at the observation location $\solution(x=1)$.
This range of solutions results in a range of parameter estimates, which flattens the posterior, but also results in a more accurate posterior mean:
because the likelihood is sharply peaked, marginalization over the meshes is approximately equal to selecting the mesh that produces the largest likelihood.
This implicitly selects the mesh that has the lowest discretization error at the observation location, and therefore results in a more accurate posterior mean than FEM.
At the same time, we see that the ground truth is not covered by the posterior distribution, due to the fact that even the best 4-element mesh will underestimate the solution at the observation location $\solution(x=1)$.
Although RM-FEM is beneficial compared to FEM for this problem, it still produces overconfident parameter estimates.
For $\meshsize=\sfrac{1}{16}$ and $\meshsize=\sfrac{1}{64}$ we see the effect of mesh randomization decrease and the RM-FEM posterior converge to the exact posterior.

Finally, we consider the statFEM posterior shown in \cref{subfig:pullout-bar-inverse-statfem} and observe no convergence of the posterior with mesh refinement.
This is problematic, but not unexpected:
a single observation is too little data to infer the two parameters of interest, let alone the three additional hyperparameters that are needed to learn the discretization error through statFEM.
As explained in the introduction to this section, it is not at all uncommon in material sciences to have to infer multiple material parameters using only force-displacement measurements at the loading point.
In such settings, in order to account consistently for discretization error, it has to be modeled rather than learned.

%% file: sections/three-point-hole.tex
\subsection{Three-point bending test}
\label{subsec:three-point-hole}
\begin{figure}
    \centering
    \includegraphics[width=\textwidth]{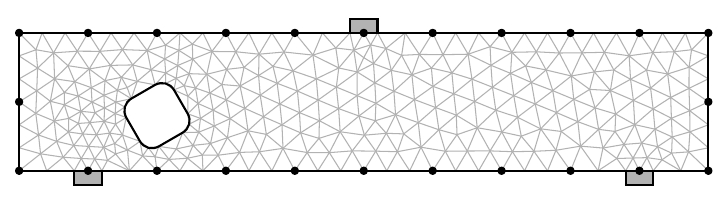}
    \caption{
        Schematic overview of the three-point bending test, with the hole in the true location.
        The coarsest mesh ($\meshsize=\qty{0.20}{\meter}$) is shown in gray.
    }
    \label{fig:three-point-hole-overview}
\end{figure}

For a 2D example, we consider the three-point bending test shown in \cref{fig:three-point-hole-overview}.
A beam with height $\beamheight=\qty{1}{\meter}$, span $\beamlength=\qty{4}{\meter}$ and overhang $\beamoverhang=\qty{0.5}{\meter}$ contains a hole at an unknown location $(\holelocx,\holelocy)$.
To determine the location of the hole, we put 24 sensors on the boundary of the beam with a uniform spacing $\observationspacing=\qty{0.5}{\meter}$.
We apply a prescribed displacement $\prescribeddisplacement=\qty{0.01}{\meter}$ at midspan and measure the resulting displacement in $x$ and $y$ direction with observation noise $\noisescale=\qty{e-4}{\meter}$.
The beam is assumed to have a constant Young's modulus $\youngsmodulus=\qty{30}{\giga\pascal}$ and Poisson's ratio $\poissonsratio=\num{0.2}$, which are typical elastic parameters of concrete.
To avoid singularities, rectangular supports of $\qty{0.20}{\meter} \times \qty{0.10}{\meter}$ with a large stiffness $\youngsmodulus_s=\qty{e6}{\giga\pascal}$ are used at the boundaries.

We assume the presence of a rounded square hole with 5 unknown parameters: center location $(\holelocx, \holelocy)$, side length $\holesize$, angle $\holeangle$ and relative corner radius $\holeradius$.
The absolute corner radius is equal to $\holeradius \holesize$, which implies that $\holeradius=0$ describes a square and $\holeradius=0.5$ describes a circle, regardless of side length $\holesize$.
The ground truth (shown in \cref{fig:three-point-hole-overview}) is fixed at $(\holelocx, \holelocy, \holesize, \holeangle, \holeradius) = (\qty{1.0}{\meter}, \qty{0.4}{\meter}, \qty{0.4}{\meter}, \sfrac{\pi}{6}, \num{0.25})$, and synthetic observations are generated by solving the forward problem with a very fine mesh ($\meshsize=\qty{0.002}{\meter}$) and corrupting the displacement at the measurement locations with \iid observation noise ($\noisescale=\qty{e-4}{m}$).
We use GMSH \cite{geuzaine_gmsh_2009} to generate finite element meshes consisting of linear triangular elements with uniform mesh size $\meshsize$.

Unlike the pullout test from \cref{subsec:pullout-bar}, this problem introduces a dependency between the parameters being inferred and the discretization itself, since a different location or shape of the hole will yield a different discretization.
In the probabilistic graph models shown in \cref{fig:graphical-models}, one should imagine an additional dependency from $\parameters$ to $\mesh$.
For RM-FEM, this additional dependency poses a problem for the Metropolis within Monte Carlo marginalization approach, since it is no longer possible to generate a set of perturbed meshes a priori and run an MCMC chain for each.
However, the Monte Carlo within Metropolis pseudomarginal method we have chosen can still be used to approximate the likelihood.
For each proposed location of the hole, $\pseudomarginalsamplecount$ perturbed meshes are generated to compute and average the likelihood.
The need to remesh also introduces discontinuities in the likelihood for all methods, since a small variation in the parameters $\parameters$ can produce a very different mesh, and therefore a jump in the likelihood.
This creates issues for samplers that assume differentiability of the target distribution, like Hamiltonian Monte Carlo, but does not pose a problem for our basic random walk Metropolis sampler.

\subsubsection{Inverse problem}
\label{subsubsec:three-point-inverse}
For the inverse problem, we put a uniform prior on each parameter:
\begin{align}
    \probabilityof*{\holelocx} &= \uniform*{0}{5} &
    \probabilityof*{\holelocy} &= \uniform*{0}{1} &
    \probabilityof*{\holesize} &= \uniform*{0}{0.5} \\
    \probabilityof*{\holeangle} &= \uniform*{0}{2 \pi} &
    \probabilityof*{\holeradius} &= \uniform*{0}{0.5} \nonumber
\end{align}
Some parameter value combinations yield problematic meshes, for example because the hole intersects a support, or cuts off a section of the beam.
To avoid this, we reject any samples where the hole intersects the boundary of the beam, which is equivalent to conditioning the prior on the fact that the hole is fully enclosed by the beam.

\begin{figure}
    \centering
    \begin{tikzpicture}
        \node (fem020) at (0, 0) {\includegraphics[width=0.25\textwidth]{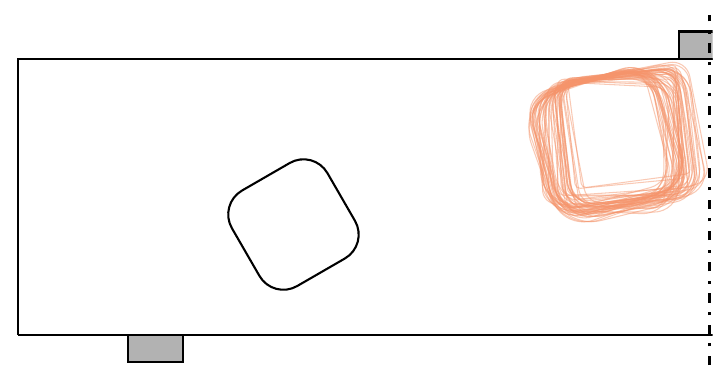}};
        \node [right=0.15 of fem020.east] (fem010) {\includegraphics[width=0.25\textwidth]{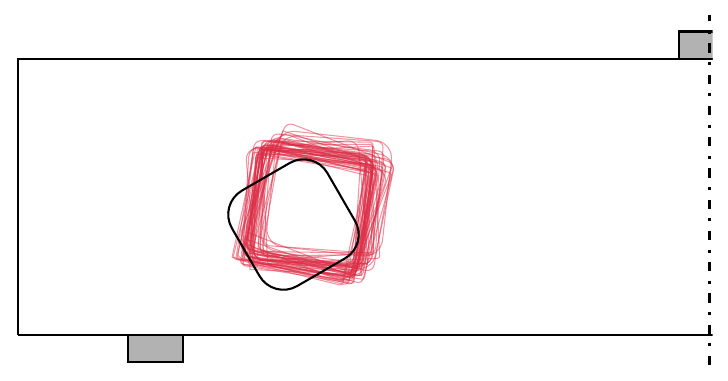}};
        \node [right=0.15 of fem010.east] (fem005) {\includegraphics[width=0.25\textwidth]{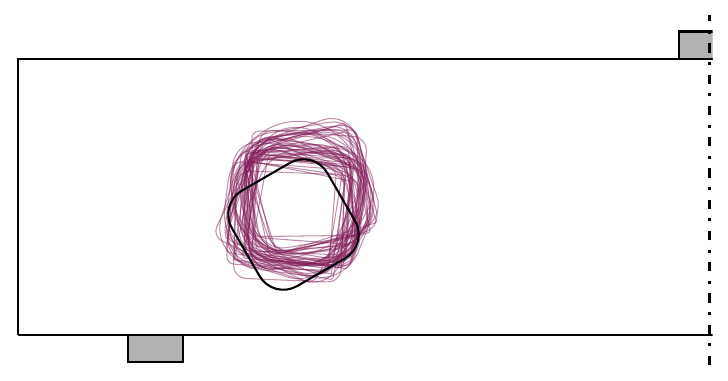}};
        \node [left=0.1 of fem020.west] {FEM};
        \node [above=0.02 of fem020.north] {$\meshsize=\qty{0.20}{\meter}$, $\elementcount\approx\num{332}$};
        \node [above=0.02 of fem010.north] {$\meshsize=\qty{0.10}{\meter}$, $\elementcount\approx\num{699}$};
        \node [above=0.02 of fem005.north] {$\meshsize=\qty{0.05}{\meter}$, $\elementcount\approx\num{2589}$};

        \node [below=0.1 of fem020.south] (bfem020) {\includegraphics[width=0.25\textwidth]{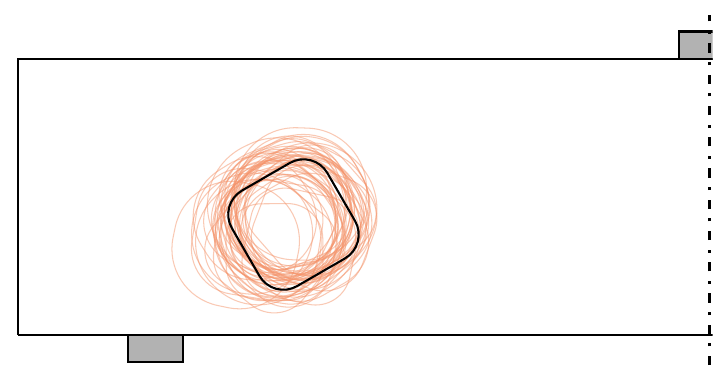}};
        \node (bfem010) at (fem010 |- bfem020) {\includegraphics[width=0.25\textwidth]{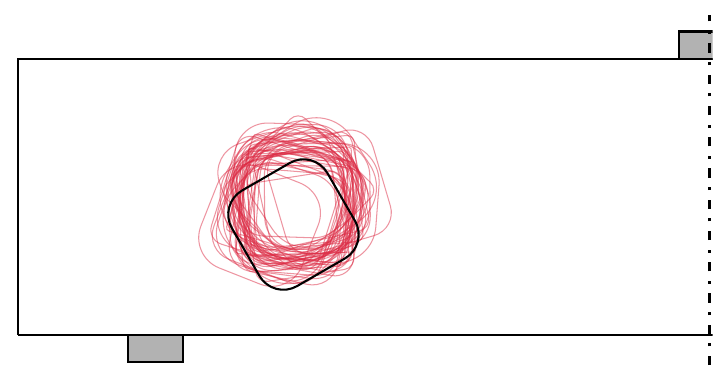}};
        \node (bfem005) at (fem005 |- bfem020) {\includegraphics[width=0.25\textwidth]{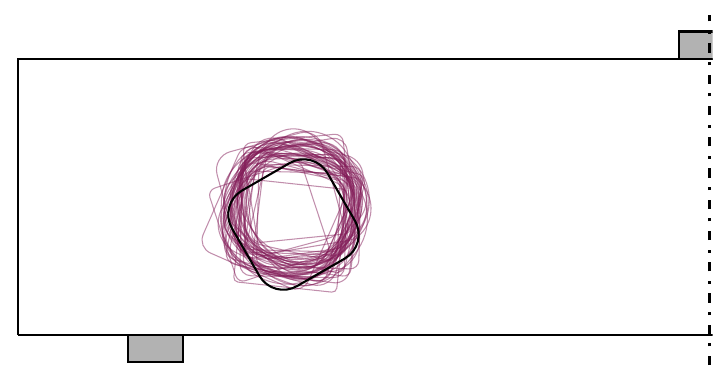}};
        \node [left=0.1 of bfem020.west] {BFEM};

        \node [below=0.1 of bfem020.south] (rmfem020) {\includegraphics[width=0.25\textwidth]{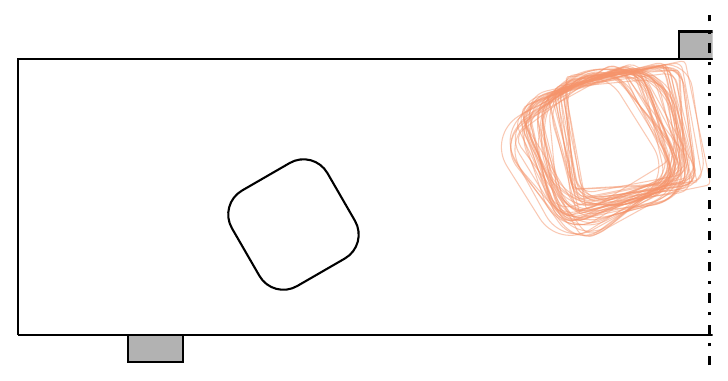}};
        \node (rmfem010) at (fem010 |- rmfem020) {\includegraphics[width=0.25\textwidth]{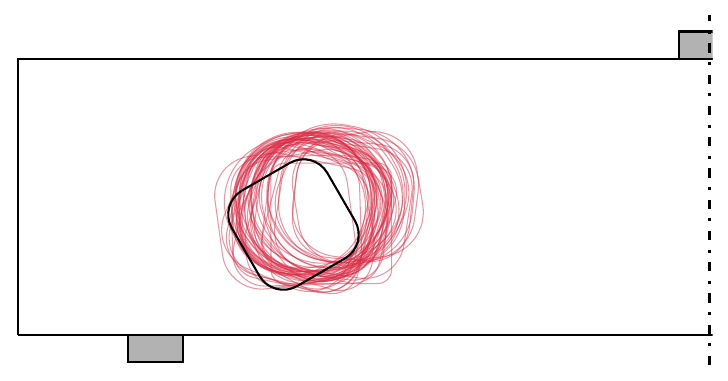}};
        \node (rmfem005) at (fem005 |- rmfem020) {\includegraphics[width=0.25\textwidth]{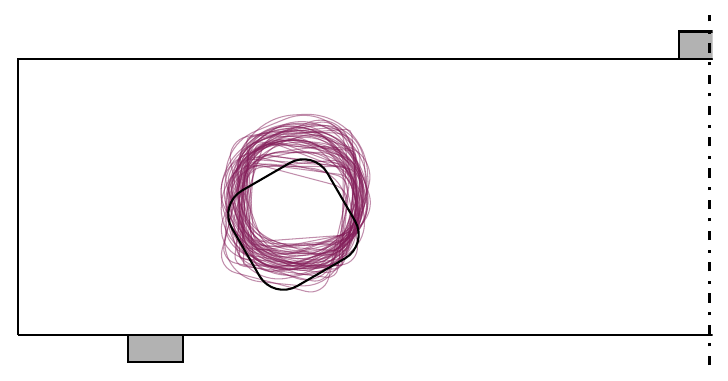}};
        \node [left=0.1 of rmfem020.west] {RM-FEM};

        \node [below=0.1 of rmfem020.south] (statfem020) {\includegraphics[width=0.25\textwidth]{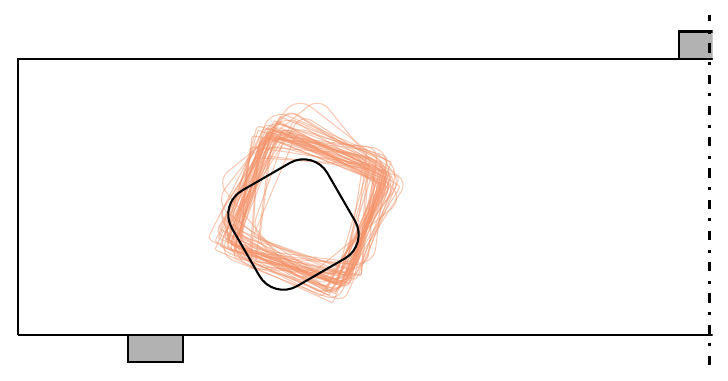}};
        \node (statfem010) at (fem010 |- statfem020) {\includegraphics[width=0.25\textwidth]{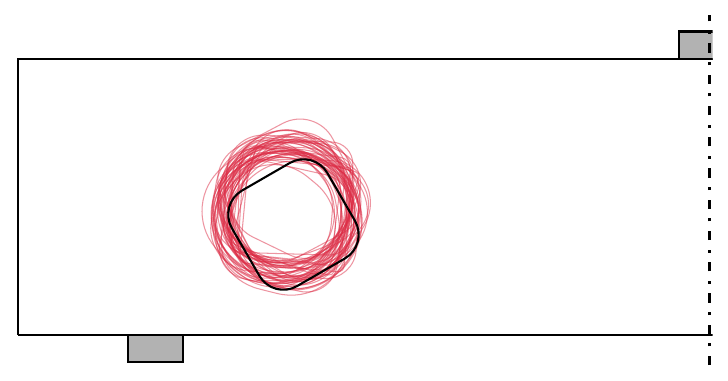}};
        \node (statfem005) at (fem005 |- statfem020) {\includegraphics[width=0.25\textwidth]{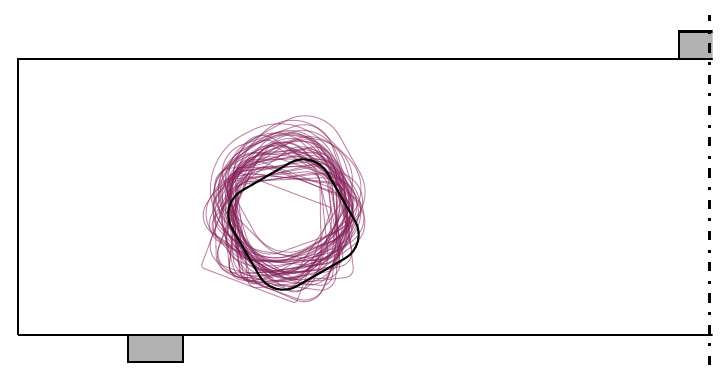}};
        \node [left=0.1 of statfem020.west] {statFEM};
    \end{tikzpicture}
    \caption{
        Samples from the posterior distribution $\posterior$ of the three-point bending test inverse problem.
        Each row corresponds to a different method, and each column corresponds to a different mesh density.
        The true location of the hole is indicated in black.
        The reported number of elements $\elementcount$ is computed with the true hole location.
    }
    \label{fig:three-point-hole-samples}
\end{figure}
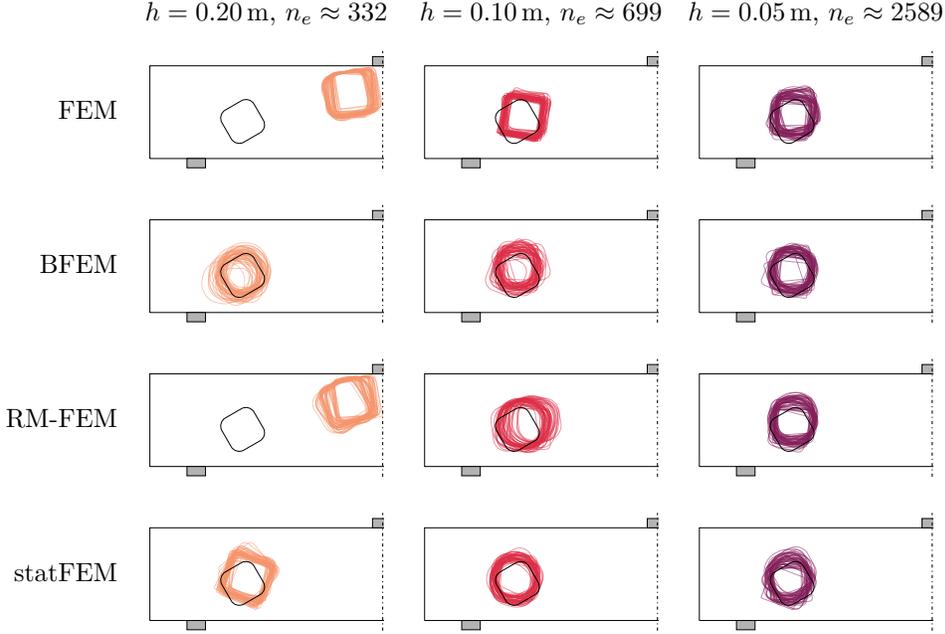

\begin{figure}
    \centering
    \includegraphics[width=\textwidth]{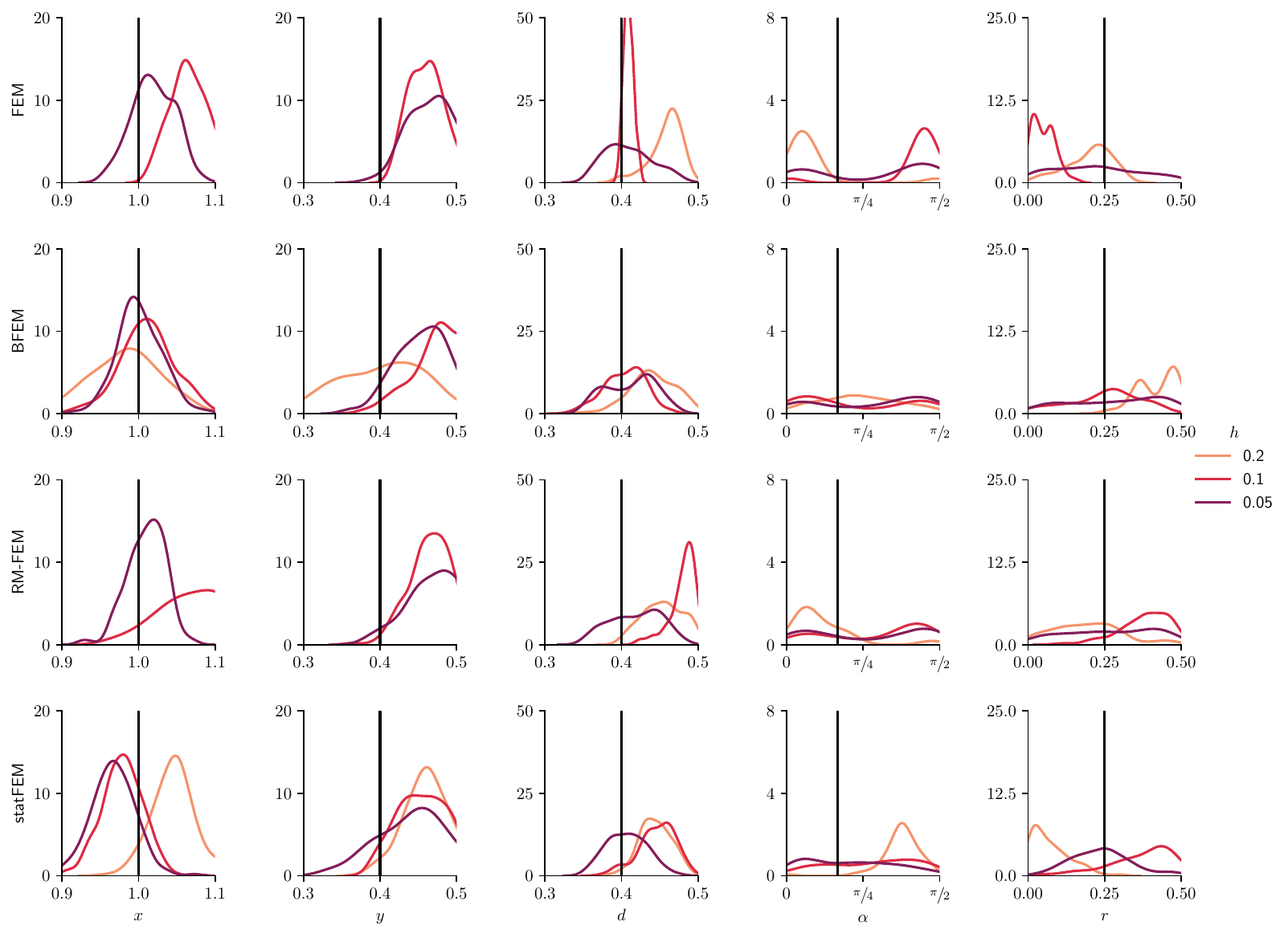}
    \caption{
        Marginal posterior distributions $\probabilityof[\observations]{\parameter_i}$ of the three-point bending test inverse problem.
        Each row corresponds to a different method, and each column corresponds to a different variable $\parameter_i$.
        The true parameter values are indicated with a black, vertical line.
    }
    \label{fig:three-point-hole-marginals}
\end{figure}

We run the inverse problem using FEM, BFEM, RM-FEM and statFEM for three different mesh sizes $\meshsize=\braces{\qty{0.20}{\meter}, \qty{0.10}{\meter}, \qty{0.05}{\meter}}$, and plot the samples of the posterior in \cref{fig:three-point-hole-samples}.
Looking at the FEM posterior samples we see that for $\meshsize=\qty{0.20}{\meter}$, there is a severe error in the inferred the location of the hole.
As the mesh is refined, the posterior mean moves towards the true location of the hole.
Using the FEM marginal distributions in the first row in \cref{fig:three-point-hole-marginals}, we can make a distinction between identifiable and unidentifiable variables.
We classify the hole location $(\holelocx, \holelocy)$ and size $\holesize$ as identifiable variables, because their posterior mean converges to the true location as the mesh is refined.
On the other hand, we classify the hole angle $\holeangle$ and corner radius $\holeradius$ as unidentifiable variables, because their posterior distributions become flat for a sufficiently fine mesh.
For both identifiable and unidentifiable variables, a coarse mesh leads to overconfidence; the difference lies in how this overconfidence is resolved with mesh refinement.
Although this distinction is problem-dependent (\eg, given more observation locations or a lower observation noise, $\holeangle$ or $\holeradius$ may become identifiable), we make the distinction nonetheless, to facilitate discussion of the results for this example.
In \cref{sec:supplement}, \cref{fig:three-point-hole-samples-supplement,fig:three-point-hole-marginals-supplement}, the FEM posterior samples and marginal distributions are shown for a broader range of meshes.
Indeed, we see that the unidentifiable parameters $\holeangle$ and $\holeradius$ converge to an approximately uniform distribution, whereas the identifiable parameters $(\holelocx, \holelocy)$ and $\holesize$ converge to a unimodal distribution which covers the ground truth.

Looking at the BFEM posterior samples in \cref{fig:three-point-hole-samples} and marginal distributions in \cref{fig:three-point-hole-marginals}, we observe that even for the coarsest mesh $\meshsize=\qty{0.20}{\meter}$, the BFEM posterior is very accurate.
The identifiable parameters $\holelocx$, $\holelocy$ and $\holesize$ are estimated correctly by the posterior mean, and BFEM produces a flat marginal distribution for the angle $\holeangle$.
Only for the relative corner radius $\holeradius$ do we see that BFEM yields an overconfident parameter estimate.
Remarkably, BFEM with a coarse observation mesh $\meshsize=\qty{0.20}{\meter}$ (which implies a reference mesh consisting of $\elementcount=4 \cdot 332=\num{1328}$ elements) performs as well as FEM with a fine mesh of $\meshsize=\qty{0.05}{\meter}$.
In \cref{sec:supplement}, \cref{fig:three-point-hole-samples-supplement}, we see that the FEM posterior for $\meshsize=\qty{0.08}{\meter}$, which has roughly the same number of elements as the reference mesh is still overconfident about the location $\holelocy$, angle $\holeangle$ and corner radius $\holeradius$ of the hole.
We see that, given a limited computational budget, it can be beneficial to spend some computational power on modeling the error, rather than use it all to maximally reduce the error.
The BFEM posteriors with finer meshes $\meshsize=\qty{0.10}{\meter}$ and $\meshsize=\qty{0.05}{\meter}$ appear very similar to the coarse BFEM posterior.
Possibly, it is not numerical error but observation information that is the limiting factor, and as a result mesh refinement does not lead to an increase in confidence in the posterior.

The RM-FEM posteriors, on the other hand, do not appear to improve much over the FEM posterior.
Although we can see the unidentifiable parameters $\holeangle$ and $\holeradius$ flatten when inspecting the marginal distributions in \cref{fig:three-point-hole-marginals}, the posterior is still very overconfident regarding the identifiable parameters $\holelocx$, $\holelocy$ and $\holesize$.
This overconfidence is apparent even more so when we consider the RM-FEM posterior samples in \cref{fig:three-point-hole-samples}.
Looking at the samples for $\meshsize=\qty{0.20}{\meter}$, we see that the RM-FEM posterior mean provides an estimate of the hole location $(\holelocx,\holelocy)$ that is just as bad as the FEM posterior mean.
Additionally, although the RM-FEM posterior samples show a bit more variance in $\holelocx$ and $\holelocy$ than the FEM posterior samples, their variance does not appear reflective of the error in the parameter estimate.

\begin{figure}
    \centering
    \begin{subfigure}{0.75\textwidth}
        \centering
        \includegraphics[width=\textwidth]{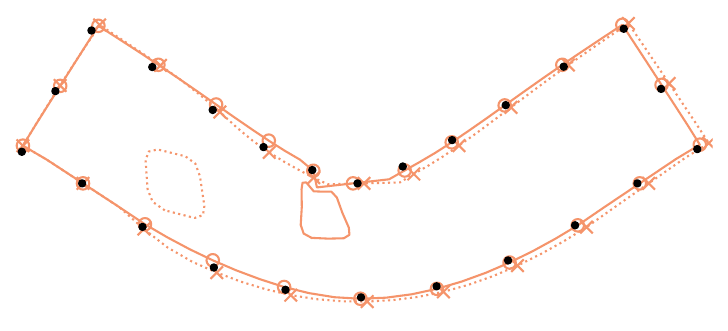}
        \caption{$h = \qty{0.20}{\meter}$}
        \label{subfig:three-point-hole-boundary-h-0.20}
    \end{subfigure}

    \begin{subfigure}{0.75\textwidth}
        \centering
        \includegraphics[width=\textwidth]{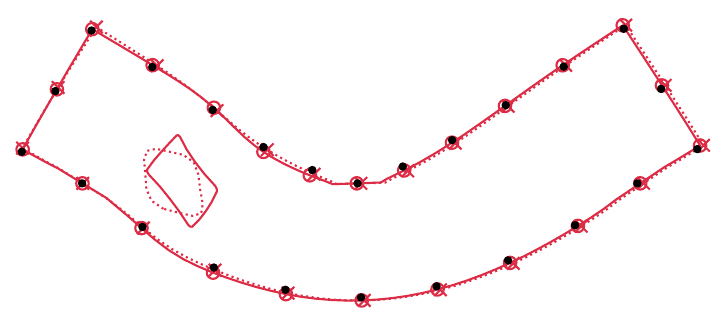}
        \caption{$h = 0.10$}
        \label{subfig:three-point-hole-boundary-h-0.10}
    \end{subfigure}

    \begin{subfigure}{0.75\textwidth}
        \centering
        \includegraphics[width=\textwidth]{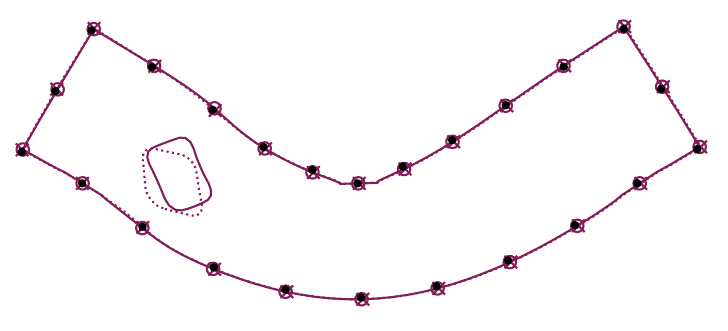}
        \caption{$h = 0.05$}
        \label{subfig:three-point-hole-boundary-h-0.05}
    \end{subfigure}
    \caption{
        FEM solution on the boundary for different mesh densities.
        Displacements are scaled with a factor $100$.
        The observed displacements $\observations$ are indicated with black dots.
        The FEM solution with the hole in the true location is indicated with a dotted line on the boundary and crosses at the observation locations.
        The FEM solution with the hole in the posterior mean point is indicated with a solid line on the boundary and circles at the observation locations.
     }
    \label{fig:three-point-hole-boundary}
\end{figure}

To better understand why RM-FEM is underperforming on this problem, we first have to consider why FEM is incorrectly inferring the location of the hole.
In \cref{fig:three-point-hole-boundary}, the predicted displacement at the boundary with the hole in the true location is compared to that of the inferred location for $\meshsize=\{\qty{0.20}{\meter}, \qty{0.10}{\meter}, \qty{0.05}{\meter}\}$.
For $\meshsize=\qty{0.20}{\meter}$, we see that the FEM prediction with the hole in the true location (crosses) does not match the observed displacements that were generated with a very fine mesh (black dots).
In particular, an increased horizontal displacement can be observed for the coarse FEM solution, which is compensated in the inverse setting by positioning the hole close to the load point.
This causes a collapse in the compressive region at the top of the beam, which reduces the horizontal motion and results in a prediction (circles) that better matches the observed displacements.
The mesh perturbation of RM-FEM does not intrinsically resolve the mechanism that causes the poor FEM posterior estimate:
each perturbed mesh produces an equally poor posterior prediction, and as a result the RM-FEM posterior appears very similar to the FEM posterior.

Finally, we consider the posterior samples of statFEM.
Looking at the marginal distributions of the unidentifiable parameters $\holeangle$ and $\holeradius$ in \cref{fig:three-point-hole-marginals}, we see that statFEM flattens the posterior compared to FEM, which resolves the overconfidence for these parameters.
At the same time, statFEM is still overconfident with regard to the identifiable parameters $\holelocx$, $\holelocy$ and $\holesize$, especially for $\meshsize=\qty{0.20}{\meter}$.
In \cref{fig:three-point-hole-samples}, we see that the statFEM posterior mean provides a better estimate of the hole location than FEM, but the ground truth still falls outside the posterior distribution.

\begin{figure}
    \centering
    \includegraphics[width=0.8\textwidth]{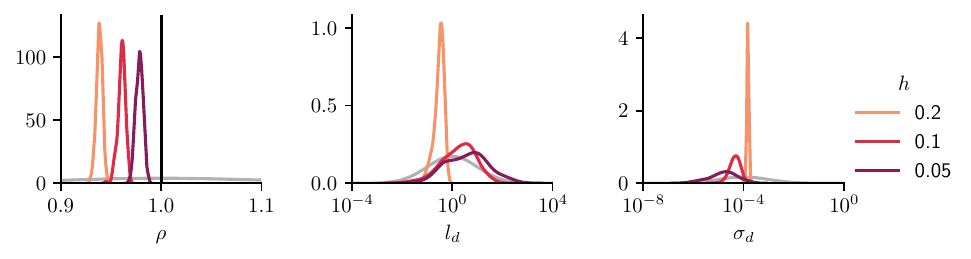}
    \caption{
        Marginal posterior distributions $\probabilityof[\observations]{\hyperparameters_i}$ of the statFEM hyperparameters for the three-point bending test inverse problem.
        The true value of the statFEM scaling parameter $\statfemscale=1.0$ is indicated with a black, vertical line.
        The prior distributions of the hyperparameters are shown in gray.
    }
    \label{fig:three-point-hole-statfem}
\end{figure}

Unlike BFEM and RM-FEM, that statFEM inference procedure involves three additional hyperparameters that were inferred from the observations, whose marginal distributions are shown in \cref{fig:three-point-hole-statfem}.
For $\meshsize=\qty{0.20}{\meter}$, we see that the statFEM scaling parameter $\statfemscale$ is estimated to be $\statfemscale \approx \num{0.933}$, which implies a reduction of the predicted displacement at the observation locations.
This reduction serves as a compensation for the overestimation of horizontal displacements by FEM as shown in \cref{fig:three-point-hole-boundary}.
As the mesh is refined, the predicted horizontal displacement is overestimated less severely, and $\statfemscale$ is estimated closer to $1$ accordingly.
Since discretization error cannot be modeled solely with a scaling factor, a misspecification component is also inferred.
For $\meshsize=\qty{0.20}{\meter}$, the length scale and magnitude of the model misspecification component $\misspecification$ are estimated to be $\misspecificationlength\approx\qty{0.333}{\meter}$ and $\misspecificationscale\approx\qty{1.74e-4}{\meter}$.
From the sharpness of the peaks for both of these hyperparameters, we can conclude that some structure of the error is being learned, with a spatial variation that is similar to the spacing of the observation locations $\observationspacing=\qty{0.50}{\meter}$ and an impact on the model prediction that is of a similar order as the observation noise $\noisescale=\qty{e-4}{\meter}$.
As the mesh is refined, this sharpness decreases, and for $\meshsize=\qty{0.05}{\meter}$, the posterior misspecification length $\misspecificationlength$ equals its prior, and the posterior misspecification scale is flat, but below $\qty{e-4}{\meter}$.
For such a fine mesh, the statFEM model is unable to learn any structure to the discretization error, but is able to tell that its impact on the prediction is less than that of the observation noise $\noise$.

%% file: sections/conclusion.tex
\section{Conclusion}
\label{sec:conclusion}
We investigated the application of BFEM to inverse problems and compared its performance to that of RM-FEM and statFEM.
In \cref{tbl:summary}, the results are summarized of the two sets of numerical experiments that have been performed:
a 1D, low-data pullout test to infer material properties, and a 2D, high-data three-point bending test to infer the geometry of a defect.
In both experiments, we observed that the discretization error of the FEM approximation to the forward model can have a significant impact on the parameter estimates that are obtained.

\begin{table}
    \centering
    \captionsetup{width=0.8\textwidth}
    \caption{
        Summary of the performance of all methods investigated in this paper.
        Overconfidence and underconfidence are checked according to \cref{subsec:misspecification-definition}.
        If a method performed differently on the pullout test and the three-point bending test, this is indicated with \badcheck/\goodcheck.
        *Convergence of statFEM depends on sufficient availability of observation data.
    }
    \label{tbl:summary}
    \begin{tabular}{lccl}
        \cmidrule[\heavyrulewidth]{2-4}
                           & not overconfident & not underconfident & convergence type \\
        \midrule
        FEM     & \badcheck & \goodcheck & FEM \textrightarrow \, exact \\
        BFEM    & \goodcheck & \goodcheck & prior \textrightarrow \, exact \\
        RM-FEM  & \badcheck & \goodcheck & FEM \textrightarrow \, exact \\
        statFEM & \goodcheck & \badcheck/\goodcheck & prior \textrightarrow \, exact* \\
        \bottomrule
    \end{tabular}
\end{table}

The BFEM approach treats discretization error in a consistent Bayesian manner:
if the FEM solution involves a severe level of error, the model considers it to provide little information regarding the inferred parameters, and defaults to the prior.
As the mesh is refined, the uncertainty surrounding the FEM solution is reduced, and the parameter estimates converge to the exact posterior.
Although BFEM is computationally expensive, because it relies on a refined mesh to define the prior distribution, it can yield a significant improvement in the accuracy of the estimated parameters.
In particular, we observed that BFEM with a coarse observation mesh and a fine reference mesh can produce a more accurate posterior than FEM with an equally fine mesh as the BFEM reference mesh.

On the other hand, RM-FEM takes a frequentist, bootstrapping approach to modeling discretization error.
The error is modeled by perturbing the node locations of the mesh, which produces a distribution over the displacement field.
The marginals of this distribution appear more similar to the error than those of BFEM, and error estimates can be derived from this distribution.
However, in the context of inverse problems, RM-FEM does not perform well, mainly due to bias in the FEM solution.
If FEM systematically over- or underestimates an observed quantity, then it will do so for each perturbed mesh.
As a result, bias is not propagated to the posterior, and RM-FEM can still produce severely overconfident parameter estimates.

Finally, statFEM is able to account for discretization error by learning it from data rather than modeling it explicitly.
Given sufficient data, the adverse effects of discretization error on the FEM prediction in the forward model can be compensated by inferring it as a form of model misspecification.
However, in low data settings, the additional hyperparameters can complicate the inference procedure or make it infeasible altogether.
In such contexts, modeling the error explicitly is more practical than inferring it in hindsight.

We have only investigated a single choice of prior for BFEM, a single perturbation scheme for RM-FEM and a single type of misspecification component for statFEM.
In this sense, the scope of the work has been quite narrow, and one could question the extent to which the results obtained for each method are fully representative of the Bayesian, frequentist and model misspecification approaches, respectively.
For example, different choices of BFEM prior exist which yield a posterior in the space orthogonal to the FEM basis, and it is an open question how they would compare to the one used in this paper.
Nonetheless, we argue that the results for each method are at least partially transferrable to their respective categories.
In particular, the RM-FEM results more broadly demonstrate the importance of modeling the error in the appropriate space, orthogonal to the FEM basis.
Likewise, the statFEM results can be extrapolated to indicate a benefit to modeling discretization error at the source where possible, rather than inferring it downstream.

Since nonlinear formulations have not yet been developed for BFEM and RM-FEM, we have limited the study to linear elliptic PDEs.
Additionally, we have left possible combinations of BFEM, RM-FEM and statFEM outside the scope of this work, to not muddle the comparison and discussion.
Combining BFEM and statFEM would be straightforward, and is potentially useful for problems where both discretization error and other sources of model misspecification play a significant role.
The mesh perturbation of RM-FEM could also be applied to BFEM, yielding multiple random observation meshes to improve the estimate of the PDE solution, or multiple random reference meshes to probe the posterior covariance operator in different trial spaces.
We leave the exploration of such ideas to future work.

%% file: sections/backmatter.tex
\appendix

\section{Finite element method for linear elasticity}
\label{sec:elasticity}
Here, we present a discussion of FEM and BFEM in the context of linear elastic solid mechanics.
On a domain $\domain \subset \reals[\dimensionality]$ with a boundary $\boundary$, we aim to find the $\dimensionality$-dimensional displacement field $\displacementfield\ofvx$.
The strain field $\strainfield$ is given by the symmetric gradient of the displacement field $\displacementfield$, the stress field $\stressfield$ is linked to the strain field $\strainfield$ via a fouth-order elasticity tensor $\elasticitytensor$, and the body force field $\sourcefield$ is given by the negative divergence of the stress field $\stressfield$:
\begin{align}
    \label{eq:elasticity-relationships}
    \strainfield &= \frac{1}{2} \parens{\gradientof \displacementfield + \T{\gradientof \displacementfield}} = \symmetricgradientof \displacementfield&
    \stressfield &= \elasticitytensor : \strainfield &
    \sourcefield &= - \divergenceof \stressfield
\end{align}
The boundary $\boundary$ is decomposed in a Dirichlet boundary $\dirichletboundary$ and a Neumann boundary $\neumannboundary$, such that $\dirichletboundary \cup \neumannboundary = \boundary$ and $\dirichletboundary \cap \neumannboundary = \emptyset$.
Combined, this gives the strong form of the PDE:
\begin{equation}
    \label{eq:elasticity-strong-form}
    \begin{aligned}
        -\divergenceof \elasticitytensor : \symmetricgradientof \displacementfield &= \sourcefield && \text{ in } \domain \\
        \displacementfield &= \dirichletfield && \text{ on } \dirichletboundary \\
        \normalfield \cdot \elasticitytensor : \symmetricgradientof \displacementfield  &= \neumannfield && \text{ on } \neumannboundary
    \end{aligned}
\end{equation}
Here, $\normalfield$ is the boundary unit normal vector, $\dirichletfield$ is the imposed displacement on the Dirichlet boundary and $\neumannfield$ is the imposed traction on the Neumann boundary.

We define the test function space $\testspace\ofOmega = \braces{\bracks{\Hspace{1}\ofOmega}^\dimensionality \, | \, \arbitraryvectorv\ofvx = 0 \, \forall \, \vx \in \dirichletboundary}$.
Multiplying both sides by a test function $\arbitraryvectorv\ofvx \in \testspace\ofOmega$ and integrating by parts yields the weak form of the PDE:
\begin{equation}
    \label{eq:elasticity-weak-form}
    \begin{aligned}
        \integrateover[\domain]{\symmetricgradientof \arbitraryvectorv : \elasticitytensor : \symmetricgradientof \displacementfield}{\vx}
        &= \integrateover[\domain]{\sourcefield \cdot \arbitraryvectorv}{\vx}
        + \integrateover[\neumannboundary]{\tractionfield \cdot \arbitraryvectorv}{\vx}
        & \forall \arbitraryvectorv \in \testspace\ofOmega
    \end{aligned}
\end{equation}
where we have made use of the fact that $\arbitraryvectorv\ofvx = \vnull$ on the Dirichlet boundary $\dirichletboundary$, and substituted the imposed traction $\tractionfield$, which takes care of the inhomogeneous Neumann boundary condition.
Note that $\gradientof \arbitraryvectorv$ has been replaced by $\symmetricgradientof \arbitraryvectorv$ without loss of generality, due to the symmetry of the stress tensor $\stressfield$.

We now define the trial function space $\trialspace\ofOmega = \braces{\bracks{\Hspace{1}\ofOmega}^\dimensionality \, | \, \arbitraryvectorv\ofvx = \dirichletfield \, \forall \, \vx \in \dirichletboundary}$ and decompose the displacement field as $\displacementfield \ofvx = \homdisplacementfield\ofvx + \inhomdisplacementfield\ofvx$, where $\homdisplacementfield \in \testspace\ofOmega$ and $\inhomdisplacementfield\ofvx = \dirichletfield \text{ on } \dirichletboundary$.
Substituting this decomposition into \cref{eq:elasticity-weak-form} yields:
\begin{equation}
    \begin{aligned}
        \label{eq:elasticity-weak-form-hom}
        \integrateover[\domain]{\symmetricgradientof \arbitraryvectorv : \elasticitytensor : \symmetricgradientof \homdisplacementfield}{\vx}
        &= \integrateover[\domain]{\sourcefield \cdot \arbitraryvectorv}{\vx}
        + \integrateover[\neumannboundary]{\tractionfield \cdot \arbitraryvectorv}{\vx}
        + \integrateover[\domain]{\symmetricgradientof \arbitraryvectorv : \elasticitytensor : \symmetricgradientof \inhomdisplacementfield}{\vx}
    \end{aligned}
\end{equation}
At this point, we have arrived at a PDE with homogeneous boundary conditions, and proceed with the usual procedure.
From the finite element mesh, we construct a set of test functions $\testvector\ofvx$ spanning $\femtestspace\ofOmega \subset \testspace\ofOmega$ and search $\femhomdisplacementfield\ofvx$ in the same space.
The entries of the resulting stiffness matrix $\stiffnessmatrix$ and force vector $\sourcevector$ are given by:
\begin{equation}
    \begin{aligned}
        \stiffnessmatrix_{ij} &= \integrateover[\domain]{\symmetricgradientof \testvector_i : \elasticitytensor : \symmetricgradientof \testvector_j}{\vx} \\
        \sourcevector_{i} &= \integrateover[\domain]{\sourcefield \cdot \testvector_i}{\vx}
        + \integrateover[\neumannboundary]{\tractionfield \cdot \testvector_i}{\vx}
        + \integrateover[\domain]{\symmetricgradientof \testvector_i : \elasticitytensor : \symmetricgradientof \inhomdisplacementfield}{\vx}
    \end{aligned}
\end{equation}
The FEM solution $\fem{\displacementfield}\ofvx$ is obtained by solving $\stiffnessmatrix \arbitraryvectorw = \sourcevector$ and reapplying the lifting function $\fem{\displacementfield}\ofvx = \arbitraryvectorw \cdot \testvector\ofvx + \inhomdisplacementfield\ofvx$.

For BFEM, we work with the homogenized weak from \cref{eq:elasticity-weak-form-hom}, and include inhomogeneous boundary conditions by making the appropriate modification to the observed force vector.
The resulting posterior distribution has a mean in $\femtestspace\ofOmega$, and a covariance that spans the subspace of $\testspace\ofOmega$ orthogonal to $\femtestspace\ofOmega$.
To obtain the posterior for the inhomogeneous problem, we simply add the lifting function $\fem{\inhomdisplacementfield}\ofvx$ to the homogeneous posterior.

\section{Supplementary results}
\label{sec:supplement}
\begin{figure}
    \centering
    \begin{tikzpicture}
        \node (fem020) at (0, 0) {\includegraphics[width=0.25\textwidth]{img/three-point-hole/sample-plot_fem_h-0.20.pdf}};
        \node [right=0.15 of fem020.east] (fem010) {\includegraphics[width=0.25\textwidth]{img/three-point-hole/sample-plot_fem_h-0.10.pdf}};
        \node [right=0.15 of fem010.east] (fem008) {\includegraphics[width=0.25\textwidth]{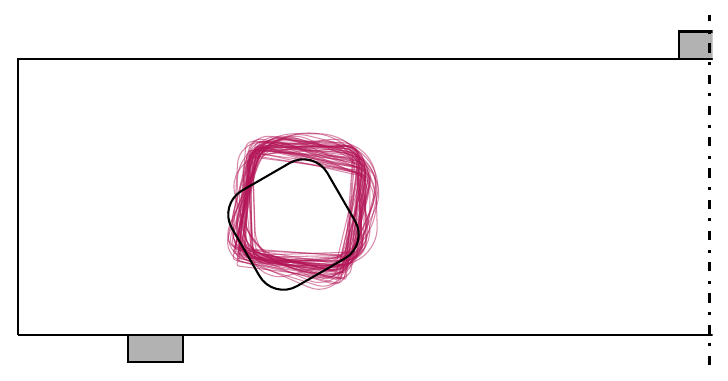}};
        \node [above=0.02 of fem020.north] {$\meshsize=\qty{0.20}{\meter}$, $\elementcount\approx\num{332}$};
        \node [above=0.02 of fem010.north] {$\meshsize=\qty{0.10}{\meter}$, $\elementcount\approx\num{699}$};
        \node [above=0.02 of fem008.north] {$\meshsize=\qty{0.08}{\meter}$, $\elementcount\approx\num{1306}$};

        \node [below=1.00 of fem020.south] (fem005) {\includegraphics[width=0.25\textwidth]{img/three-point-hole/sample-plot_fem_h-0.05.pdf}};
        \node (fem002) at (fem010 |- fem005) {\includegraphics[width=0.25\textwidth]{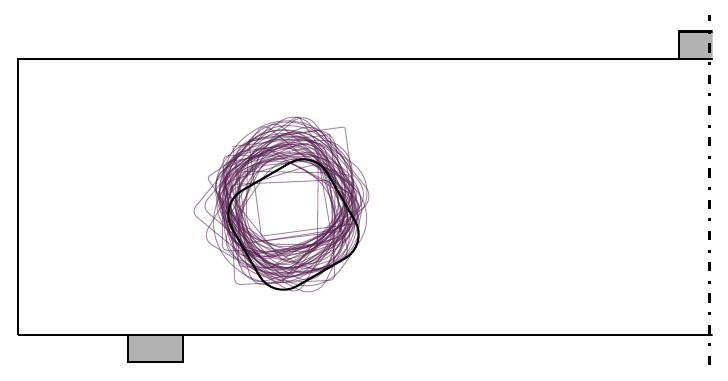}};
        \node (fem001) at (fem008 |- fem005) {\includegraphics[width=0.25\textwidth]{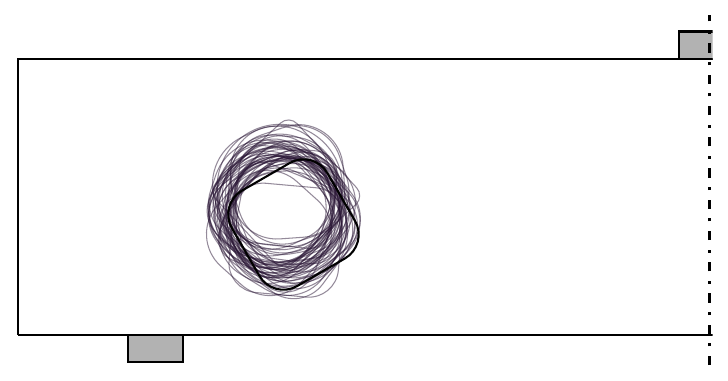}};
        \node [above=0.02 of fem005.north] {$\meshsize=\qty{0.05}{\meter}$, $\elementcount\approx\num{2589}$};
        \node [above=0.02 of fem002.north] {$\meshsize=\qty{0.02}{\meter}$, $\elementcount\approx\num{14720}$};
        \node [above=0.02 of fem001.north] {$\meshsize=\qty{0.01}{\meter}$, $\elementcount\approx\num{57933}$};
    \end{tikzpicture}
    \caption{
        Samples from the FEM posterior distribution $\posterior$ of the three-point bending test inverse problem for six different mesh sizes $\meshsize$.
        The true location of the hole is indicated in black.
    }
    \label{fig:three-point-hole-samples-supplement}
\end{figure}

\begin{figure}
    \centering
    \includegraphics[width=\textwidth]{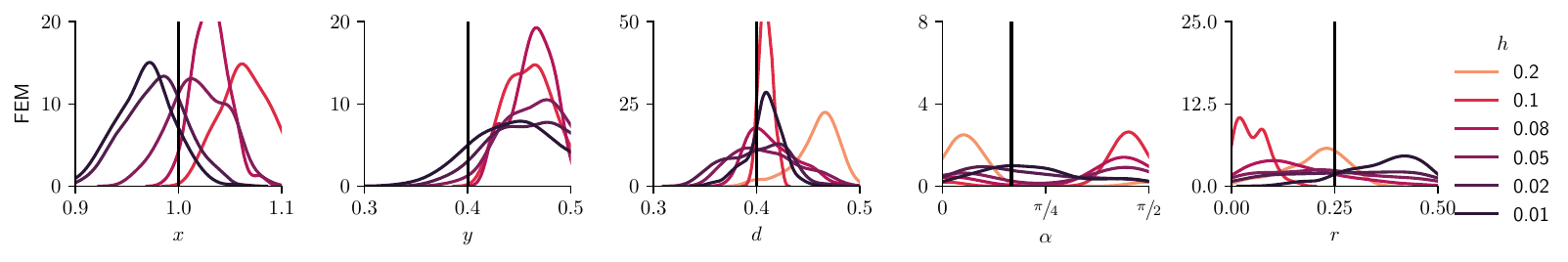}
    \caption{
        Marginal FEM posterior distributions $\probabilityof[\observations]{\parameter_i}$ of the three-point bending test inverse problem for six different mesh sizes $\meshsize$.
        The true parameter values are indicated with a black, vertical line.
    }
    \label{fig:three-point-hole-marginals-supplement}
\end{figure}

To complement the results from \cref{subsec:three-point-hole}, we present the FEM posteriors for a broader range of mesh densities.
In \cref{fig:three-point-hole-samples-supplement}, the FEM posterior samples are shown for six different mesh sizes, similar to \cref{fig:three-point-hole-samples}
In \cref{fig:three-point-hole-marginals-supplement}, the corresponding posterior marginal distributions are shown, similar to \cref{fig:three-point-hole-marginals}.

\bibliographystyle{bst/siamplain}
\bibliography{citations}